\documentclass[11pt]{article}

\usepackage{amsmath, amssymb, amscd, amsthm, amsfonts}
\usepackage{graphicx} 
\usepackage{hyperref} 
\usepackage{booktabs}
\usepackage{ dsfont, bm}
\usepackage[title]{appendix}
\usepackage{comment}

\usepackage[textwidth=500pt,textheight=650pt,centering]{geometry} 

\title{Continuous-time weakly self-avoiding walk on $\mathbb{Z}$ has strictly monotone escape speed}
\author{Yucheng Liu\thanks{Department of Mathematics, University of British Columbia, Vancouver, BC, Canada V6T 1Z2. 
\url{https://orcid.org/0000-0002-1917-8330}, 
\href{mailto:yliu135@math.ubc.ca}{yliu135@math.ubc.ca}.
}}
\date{\vspace{-5ex}} 

\newtheorem{theorem}{Theorem}[section]
\newtheorem{lemma}[theorem]{Lemma}

\newtheorem{proposition}[theorem]{Proposition}
\newtheorem{corollary}[theorem]{Corollary}

\theoremstyle{definition}
\newtheorem{remark}[theorem]{Remark}

\newcommand{\Z}{\mathbb{Z}}

\newcommand{\R}{\mathbb{R}}
\newcommand{\C}{\mathbb{C}}

\numberwithin{equation}{section}

\newcommand{\ie}{\textit{i}.\textit{e}.}
\newcommand{\eg}{\textit{e}.\textit{g}.}
\newcommand{\eps}{\varepsilon}

\newcommand{\N}{\mathbb{N}}
\newcommand{\E}{\mathbb{E}}

\renewcommand{\P}{\mathbb{P}}
\newcommand{\del}{\partial}
\newcommand{\inv}{^{-1}}
\renewcommand{\(}{\left(}
\renewcommand{\)}{\right)}
\newcommand{\half}{\frac{1}{2}}
\newcommand{\Del}{\Delta}
\newcommand{\Lcal}{\mathcal{L}}

\newcommand{\1}{\mathds{1}}
\newcommand{\stle}{\preccurlyeq}
\renewcommand{\bf}{\bm}

\renewcommand{\Re}{\mathrm{Re}\,}
\providecommand{\abs}[1]{\lvert#1\rvert} 
\providecommand{\norm}[1]{\lVert#1\rVert}

\newcommand{\dt}{{dt}}

\begin{document}
\maketitle

\begin{abstract}
Weakly self-avoiding walk (WSAW) is a model of simple random walk paths that penalizes self-intersections. On $\mathbb{Z}$, 
Greven and den Hollander proved in 1993 that the discrete-time weakly self-avoiding walk has an asymptotically deterministic escape speed, and they conjectured that this speed should be strictly increasing in the repelling strength parameter. 
We study a continuous-time version of the model,
give a different existence proof for the speed, and prove the speed to be strictly increasing. 
The proof uses a transfer matrix method implemented via a supersymmetric version of the BFS--Dynkin isomorphism theorem, spectral theory, Tauberian theory, and stochastic dominance. 
\end{abstract}

\section{Introduction}
Weakly self-avoiding walk is a model based on the simple random walk that penalizes self-intersections. In the discrete-time setting, it is also known as the self-repellent walk and as the Domb-Joyce model for soft polymers \cite{DombJoyce1972}. In the model, every self-intersection of the walk is penalized using a factor $1-\alpha$, with parameter $\alpha\in (0,1)$. 
The boundary value $\alpha = 0$ corresponds to the simple random walk and $\alpha=1$ corresponds to the strictly self-avoiding walk. 

In dimension one, the strictly self-avoiding case $\alpha=1$ has only two walks, going left or right, so it has escape speed equal to 1. Greven and den Hollander \cite{GH1993} proved in 1993 that the weakly self-avoiding walk also escapes linearly from the origin as length goes to infinity, with an asymptotically deterministic speed $\theta^*(\alpha)\in (0, 1)$ satisfying $\lim_{\alpha\to 0}\theta^*(\alpha) = 0$ and $\lim_{\alpha\to 1}\theta^*(\alpha) =1$. 
The result was extended to a central limit theorem by K\"onig in 1996 \cite{Konig1996}. Their proofs use large deviation theory and local times of the simple random walk. In 2001, van der Hofstad \cite{Hofstad2001} gave another proof of the central limit theorem for $\alpha$ close to 1, using the lace expansion. 

The escape speed $\theta^*(\alpha)$ is conjectured to be strictly increasing in the repelling strength parameter $\alpha$. The walk should escape faster if the self-repellency is stronger. But to our knowledge, this has not been proved. 	
In this paper, we study a continuous-time version of the model. We prove this model also has an asymptotically deterministic speed and the speed is strictly increasing in the repelling strength parameter. 
Our proof of the existence of escape speed uses local times of the simple random walk and Tauberian theory. The monotonicity of the speed is proved using stochastic dominance. 
The speed is qualitatively similar to the speed in the discrete-time model \cite{GH1993}, in the sense that they are both identified as the reciprocal of the derivative of the largest eigenvalue of some operator.

The continuous-time weakly self-avoiding walk was studied in \cite{BBS2015-0, BBS2015} on $\Z^4$, and a generalization of the model was studied in \cite{BS2020} on the complete graph,
both using supersymmetric representation. Supersymmetric representation is a way of expressing certain functionals of local times as an integral of differential forms. We give a brief introduction to the representation in Appendix~\ref{appendix:supersymmetry}. It allows us to write the two-point function of the WSAW in a nice form to apply the transfer matrix approach. 

There is an analogous continuous-time-and-space model based on Brownian motion called the Edwards' model. In this model, there is also an escape speed, first proved by Westwater in 1985 \cite{Westwater1985} using Brownian local times. The speed is very simple due to the Brownian scaling property, and it is strictly increasing in the repelling strength parameter.

\subsection{The continuous-time WSAW model and main results}
We consider the continuous-time nearest-neighbor simple random walk on $\Z$, \ie, the walk jumps left or right, each with rate 1. Denote the walk by $\{X(t)\}_{t\ge0}$ and denote its expectation by $E_i$ if $X(0)=i$. 
The \emph{local time of $X$ at $x\in \Z$ up to time $T$} is defined by \begin{equation} \label{def:local-time}
L_{T,x} = \int_0^T \1_{X(s)=x}\ ds. \end{equation}
Notice that \begin{align}
L_{T,x}^2 =  \int_0^T \int_0^T \1_{X(t)=x}\1_{X(s)=x}\ dsdt
= \int_0^T \int_0^T \1_{X(t)=X(s)=x}\ dsdt
\end{align} gives the \emph{self-intersection local time at $x$ up to time $T$}. 

As in the discrete-time setting, we penalize self-intersections. 
Let $g>0$ be the repelling strength parameter.
The weakly self-avoiding walk measure $\P_i^{g, T}$, for walks 
starting from $i$, is defined by the $\phi(t) = t^2$ case of the
expectation \begin{align}
\E_i^{g,T}[ f(X) ]\propto E_i \( e^{-g \sum_{x=-\infty}^\infty \phi(L_{T,x})} f(X) \). 
\end{align}
For $i, j\in \Z$, we define \begin{equation}
P_{ij}(g, T) = E_i \(  e^{-g\sum_{x=-\infty}^\infty\phi(L_{T,x})} \1_{X(T)=j} \).  \end{equation}
We use the Laplace transform with a complex parameter $\nu\in\C$.
With $p(t) = e^{-g\phi(t) - \nu t}$,
the \emph{two-point function} of the weakly self-avoiding walk from $i$ to $j$ is defined to be 
\begin{align}
G_{ij}(g, \nu) &= \int_0^\infty P_{ij}(g, T) e^{-\nu T}dT 
= \int_0^\infty E_i \( \prod_{x=-\infty}^\infty p(L_{T,x}) \1_{X(T)=j} \) dT,
\end{align}
where the second equality follows from $T=\sum_{x=-\infty}^\infty L_{T,x} $. 
A positive $\nu > 0$ penalizes the length of the walk and acts as a killing rate. 
The use of complex $\nu$ is mainly to apply our complex Tauberian theorem, 
which transfers asymptotics of $G_{ij}$ to asymptotics of $P_{ij}$. 
The relevant range of parameters will have $\Re(\nu) \le 0$. 
Our methods allow us to consider a more general model. We only assume the function $\phi:[0, \infty)\to [0, \infty)$ satisfies \begin{align}
&\phi(0)=0, \label{A0}\tag{H0} \\
&\phi(t) \text{ is differentiable for all } t\ge 0, \label{A1} \tag{H1} \\
&\phi(t)/t \text{ is increasing,} \label{A2} \tag{H2} \\
&\lim_{t\to\infty} \phi(t)/t  = \infty,\label{A3} \tag{H3} \\
&\phi'(t) e^{-g\phi(t) -\nu t} \text{ is a Schwartz function for all } g>0, \nu\in \R. \label{A4} \tag{H4}
\end{align}
For example, $\phi(t)$ can be a polynomial $\phi(t) = \sum_{k=2}^M a_k t^k$ with $M\ge2$, $a_M>0$, and $a_k \ge 0$ for all $k$.

By translation invariance, we can always start the walk at 0. 
The main result of the paper is the following theorem, which asserts that the weakly self-avoiding walk has an escape speed and the speed is strictly increasing in the repelling strength $g$. The theorem is stated for walks going to the right, but walks going to the left have the same speed by symmetry. 

\begin{theorem} \label{theorem:WLLN}
There exists an analytic function $\theta: (0, \infty) \to (0, \infty)$ with $\theta' > 0$ everywhere such that for all $g>0$, $\eps>0$, 
\begin{align} \label{eq:WLLN}
\lim_{T\to \infty} \P_0^{g, T} \(  \left\lvert \frac{X(T)}{T} - \theta(g) \right\rvert \ge \eps 
\,\middle\vert\, 
X(T) > 0 \) = 0. 
\end{align}
\end{theorem}

We also have the following result on the mean end-to-end distance, 
which readily implies $L^p$ convergence of $\frac{X(T)}T$ to $\theta(g)$ for any positive even integer $p$ and the convergence in probability \eqref{eq:WLLN}.
The notation $f(T)\sim h(T)$ means $\lim_{T\to\infty} f(T)/g(T)=1$. 
\begin{theorem} \label{theorem:moments}
For the same $\theta(\cdot)$ as in Theorem~\ref{theorem:WLLN}, for any $g>0$ and any $k\in \N$, \begin{align} \label{eqn:6}
\E_0^{g,T}\big[ X(T)^k \mid X(T) > 0 \big] =
\frac{ \sum_{j=1}^\infty j^k P_{0j}(g, T) }{ \sum_{j=1}^\infty P_{0j}(g, T) } \sim \( \theta(g) T\)^k, \qquad T\to \infty. \end{align}
\end{theorem}

\begin{figure}
\centering
\includegraphics[scale=0.8]{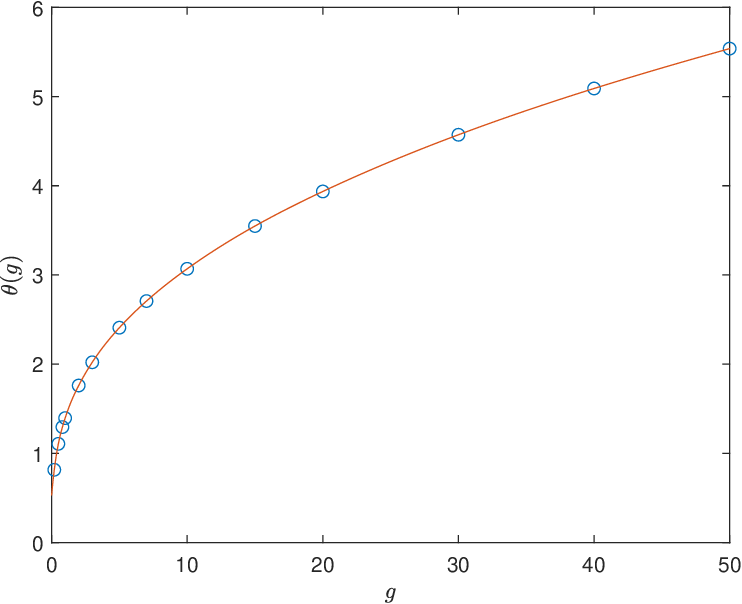}
\caption{
Numerical evaluation of the escape speed $\theta(g)$ of the WSAW model ($\phi(t) = t^2$). 
The dots are the numerical results, interpolated by cubic spline.
The computation is based on truncating and discretizing the integral operator $Q$ (defined in~\eqref{def:Q}) at $s=100$ and with step size $0.001$.  
The numerical results suggest $\theta(g) \sim Cg^{1/3}$ as $g\to 0$. 
For the discrete-time model, this asymptotic relation is proved by van der Hofstad and den Hollander in \cite{HH1995}.
}
\label{figure1}
\end{figure}
As a by-product of the proofs, we obtain critical exponents for the two-point function, the susceptibility, and the correlation lengths of integer orders. 
The critical exponents are the same as for the strictly self-avoiding walk, 
and they are collected in the following theorem. 
For all $g>0$, we will define $\nu_c(g)$ to be the $\nu \in \R$ at which an integral operator has operator norm exactly $1$, see Section \ref{section:Qproperty}. We will show this is consistent with the usual definition of $\nu_c(g)$, which is by the property that \begin{align} \label{usual-nu_c}
\chi(g, \nu) < \infty 
\quad \text{if and only if} \quad
\nu> \nu_c(g), \end{align}
where $\chi(g, \nu) = \sum_{j=-\infty}^\infty G_{0j}(g, \nu)$ is the \emph{susceptibility}.

\begin{theorem} \label{theorem:critical_exponent}
Let $g>0$, $\nu\in \R$, and $\nu_c = \nu_c(g)$ (defined in Section \ref{section:Qproperty}). 
There exist non-zero constants $A, B, C_1, C_2, \dots$ such that: 

\begin{enumerate}
\item[(i)]
The critical two-point function $G_{ij}(g, \nu_c)$ satisfies, as $\abs{ j-i }\to \infty$, \begin{align} \label{exponent:eta}
G_{ij}(g, \nu_c) \sim \frac{ A}{ \abs{ j-i }^{d-2+\eta}} = A, \qquad d=\eta=1.
\end{align}

\item[(ii)]
The susceptibility $\chi(g, \nu)$ satisfies, as $\nu\to \nu_c^+$, 
\begin{align}
\chi(g, \nu)
\sim  B (\nu-\nu_c)^{-\gamma}, \qquad \gamma=1.
\end{align}

\item[(iii)]
For any $k\in \N$, the correlation length of order $k$ satisfies, as $\nu\to \nu_c^+$, \begin{align} \label{exponent:nu_k}
\xi_k(g, \nu) 
= \(\frac{\sum_{j=-\infty}^\infty \abs j^k G_{0j}(g, \nu)}{\chi(g, \nu)}\)^{1/k}
\sim 
C_k (\nu-\nu_c)^{-\nu_k}, \qquad \nu_k = 1.
\end{align}
\end{enumerate}
\end{theorem}

The rest of the paper is organized as follows. 
In Section~\ref{section:2pt}, we first make a finite-volume approximation and use the supersymmetric representation of random walks to express the finite-volume two-point function as an inner product (Proposition~\ref{prop:Gij^N}). Then we take the infinite-volume limit and study its asymptotic behavior near the critical $\nu_c$. 
We also calculate the susceptibility and correlation lengths, and prove Theorem~\ref{theorem:critical_exponent} in Corollary~\ref{corollary:Gij}, Proposition~\ref{prop:chi}, and Corollary~\ref{corollary:correlation}. 
In Section~\ref{section:laplace}, we first prove a general Tauberian theorem, then we use the theorem, with asymptotics in the parameter $\nu$ from Section~\ref{section:2pt} as input, to prove Theorem~\ref{theorem:moments}. 
The convergence part of Theorem~\ref{theorem:WLLN} follows. 
In Section~\ref{section:monotonicity}, we use a separate stochastic dominance argument to prove $\theta' > 0$ everywhere, finishing the proof of Theorem~\ref{theorem:WLLN}.

\section{Two-point function, susceptibility and correlation length} \label{section:2pt}
In this section, we first work on the finite subset $[-N, N]\cap \Z$ with free boundary conditions. We use the transfer matrix approach to derive an expression for the finite-volume two-point function $G_{ij}^N(g, \nu)$. Then we define a critical parameter $\nu_c(g) \in \R$, and we show that for $\Re(\nu) \ge \nu_c(g)$,
the infinite-volume limit $\lim_{N\to\infty} G_{ij}^N(g, \nu)$ exists and equals $G_{ij}(g, \nu)$. 
We use this representation of $G_{ij}$ to study one-sided versions of the susceptibility and correlation lengths, obtaining asymptotics as $\nu \to \nu_c$.
Two-sided versions of the quantities follow readily by symmetry. 

\subsection{Finite-volume two-point function} \label{section:finite_volume}
Let $N\in \N$. Consider the nearest-neighbor simple random walk on $[-N, N]\cap \Z$. That is, the walk jumps to left or right if possible, each with rate $1$. For $i,j\in [-N, N]$, let $E^N_i$ denote the expectation of such a walk starting at $i$. The local times of this walk are defined exactly as in \eqref{def:local-time}. Recall $p(t) = e^{-g\phi(t) - \nu t}$. 
We define \begin{align}
P_{ij}^N (g, T) &= E_i^N \( e^{-g \sum_{x=-N}^N \phi(L_{T,x})} \1_{X(T)=j} \), \\
G_{ij}^N (g, \nu) &= \int_0^\infty P_{ij}^N (g, T) e^{-\nu T}dT 
= \int_0^\infty E_i^N \( \prod_{x=-N}^{N} p(L_{T,x}) \1_{X(T)=j} \) dT. \label{def:Gij^N}
\end{align}
The modified Bessel functions of the first kind of orders 0 and 1 are the entire functions \begin{align}
I_0(z) &= \sum_{m=0}^\infty \frac{1}{ m! m!} \left(\frac z 2\right)^{2m}, \qquad
I_1(z) = \sum_{m=0}^\infty \frac{1}{ m! (m+1)!} \left(\frac z 2\right)^{2m+1} \label{def:I0I1}
\end{align} respectively. For $j=0,1$, we define the $g, \nu$-dependent
symmetric kernels \begin{align} \label{def:k_j}
k_j(t,s)= \sqrt{p(t)}\sqrt{p(s)} e^{-t}e^{-s} I_j(2\sqrt{st}), \end{align}
where $\sqrt{p(t)} = e^{-\half g\phi(t) -\half \nu t}$.
Since $p(\cdot)$ has exponential decay by assumption~\eqref{A3}, the kernels are square-integrable. 
We also define the $g, \nu$-dependent operators $T,Q: L^2[0,\infty) \to L^2[0,\infty)$, by \begin{align}
Tf(t) 
&= \sqrt{p(t)} e^{-t} + \int_0^\infty f(s) k_1(t,s) \sqrt{\tfrac{t}{s}} ds, \label{def:T} \\
Qf(t) 
&= \int_0^\infty f(s)k_0(t,s) ds. \label{def:Q}
\end{align}
It follows from the definition of $k_1$ and the Taylor series of $I_1$ that \begin{align} \label{eqn:2.6}
k_1(t,s) \sqrt{\tfrac{t}{s}}
=  \sqrt{p(t)}\sqrt{p(s)} e^{-t}e^{-s} 
\sum_{m=0}^\infty \frac{s^m t^{m+1}}{ m! (m+1)!},
\end{align}
which is regular at $s=0$ for all $t$. Notice $T$ is affine but non-linear. 
For $Q$, since $k_0(t,s)$ is square-integrable, $Q$ is a Hilbert--Schmidt operator, thus compact. 
When $\nu\in \R$, $k_0(t,s)$ is real and symmetric, so $Q$ is self-adjoint.
With these operators, we prove the following representation of $G_{ij}^N$. The inner product here is the usual $\langle f, h \rangle = \int_0^\infty f(t)\overline{ h(t) } dt$. 
\begin{proposition} \label{prop:Gij^N}
Let $g>0$, $\nu\in \C$, $N\in \N$, and $-N\le i\le j\le N$. Then \begin{align}
G_{ij}^N(g, \nu) 
&=\left \langle Q^{j-i}T^{N+i}[\sqrt p] , \overline{ \Big( T^{N-j}[\sqrt p] \Big)}\right \rangle. 
\end{align}
\end{proposition}
Note the complex conjugation on $ T^{N-j}[\sqrt p]$ cancels with the conjugation from the inner product. 
The proposition is proved using the transfer matrix approach, implemented via a supersymmetric representation of the random walk. This is the only place we need the supersymmetric representation. The proof of the proposition and an introduction to the supersymmetric theory are included in Appendix~\ref{appendix:supersymmetry}. 

\begin{remark}
Fix $x \in \N$. For a continuous-time simple random walk starting at 0 and stopped when the local time $L_{T,x}$ reaches some fixed value, the local times between 0 and $x$ form a Markov chain with the transition kernel $e^{-t}e^{-s} I_0(2\sqrt{st})$, by the Ray--Knight Theorem \cite[Theorem 4.1]{BHK2007}. The kernel $k_0(t,s)$ for our operator $Q$ contains this Markov kernel and has the extra $\sqrt{p(t)}\sqrt{p(s)}$ to incorporate self-interactions. 
Local times outside the interval $[0,x]$ also form Markov chains, and the Markov kernel there is similarly related to the operator $T$. 
\end{remark}

\subsection{Critical parameter $\nu_c$ and properties of $Q$} \label{section:Qproperty}
In this section, we prove properties of the operator $Q$ and define the critical parameter $\nu_c(g) \in \R$. 

\begin{lemma} \label{lemma:Q}
For $g > 0$ and $\nu\in\R$:
\begin{enumerate}
\item[(i)]
The operator $Q$ defined by \eqref{def:Q} is positive, \ie, $\langle Qf,f \rangle\ge0$ for all $f$. Hence, all eigenvalues of $Q$ are non-negative.  

\item[(ii)]
The operator norm $\|Q\|$ is a simple, isolated eigenvalue of $Q$, and there exists an eigenvector with eigenvalue $\|Q\|$ that is continuous and strictly positive everywhere. 
\end{enumerate}
\end{lemma}

\begin{proof}
(i) We use the Taylor series $I_0(z) = \sum_{m=0}^\infty \frac{1}{ m! m!} (\frac z 2)^{2m}$.
For any $f$, we have \begin{align}
\langle Qf,f \rangle
&= \int_0^\infty \int_0^\infty  \overline{f(t)} f(s)  \sqrt{p(t)}\sqrt{p(s)} e^{-t}e^{-s} I_0(2\sqrt{st})\ dsdt  \\
&=  \sum_{m=0}^\infty \int_0^\infty  \overline{ f(t)} \sqrt{p(t)} e^{-t}  \frac{t^m}{m!}\ dt
\int_0^\infty   {f(s)} \sqrt{p(s)} e^{-s}  \frac{s^m}{m!}\ ds\nonumber \\
&=  \sum_{m=0}^\infty   \left\lvert\int_0^\infty  f(t) \sqrt{p(t)} e^{-t}  \frac{t^m}{m!}\ dt \right\lvert^2 
\ge 0, \nonumber
\end{align} because $\sqrt{p}$ is real for real $\nu$. 

\smallskip \noindent
(ii) Since $Q$ is compact and self-adjoint, one of $\pm\|Q\|$ is an eigenvalue \cite[Lemma 6.5]{SteinRealAnalysis}. Since eigenvalues must be non-negative by part (i), $\norm Q$ must be an eigenvalue. It is isolated by compactness. 
Also, since any eigenvalue of $Q$ must have magnitude less than the operator norm, the spectral radius $r(Q)$ equals $\norm Q$. 

We use a Krein--Rutman theorem for irreducible operators on Banach lattices \cite[Theorem~1.5]{chang2020spectral} (also see references therein). 
Krein--Rutman theorems are generalizations of the Perron--Frobenius theorem for entrywise positive matrices. 
Since $X = L^2[0, \infty)$ is a Banach lattice equipped with the cone $P$ of non-negative functions, once we verify that $Q(P) \subset P$ and that $Q$ is ideal-irreducible and compact, 
by \cite[Theorem~1.5]{chang2020spectral} we will learn that $r(Q) = \norm Q$ is a \emph{simple} eigenvalue, with an eigenvector $h$ that is a quasi-interior point of the cone. 
By (ii) in the third paragraph before \cite[Theorem~1.5]{chang2020spectral}, 
$h$ is a quasi-interior point of $P$ if and only if $\langle f, h \rangle > 0$ for all nonzero $f \ge 0$, which happens if and only if $h$ is strictly positive everywhere. 
Moreover, since
\begin{align}
\|Q\| h(t) = (Qh)(t) = \int_0^\infty h(s) k_0(t,s)ds
\end{align}
and $k_0$ is continuous, $h$ is continuous too.

It remains to verify the conditions on $Q$. $Q$ is compact because it is a Hilbert--Schmidt operator, and $Q(P) \subset P$ because $k_0\ge 0$.
To check that $Q$ is ideal-irreducible, we need to show that any nonzero $Q$-invariant closed ideal\footnote{An ideal $I \subset X$ is a linear subspace of $X$ with the property that $y\in I$, $x\in X$, and $\abs x \le \abs y$ imply $x\in I$.}
$I \subset X$ equals $X$. 
Let $0 \ne u \in I$. By the definition of ideal, we have $\abs u \in I$. 
By $Q$-invariance, $Q(\abs u) \in I$ also. 
Since $k_0 > 0$ and $\abs u \ne 0$, we have $Q(\abs u) >0$ everywhere, which makes it a quasi-interior point of $P$. 
Since the ideal generated by a quasi-interior point is dense in $X$ by the third paragraph before \cite[Theorem~1.5]{chang2020spectral}, we must have $I = X$, as $I$ is closed. 
This establishes ideal-irreducibility and concludes the proof. 
\end{proof}

We define $\nu_c(g)$ to be the $\nu\in \R$ that satisfies $\|Q(g, \nu)\| =1$. The following lemma guarantees $\nu_c(g)$ exists and is unique. 
In Proposition~\ref{prop:chi} and the paragraphs after it, we will prove $\chi(g, \nu)<\infty$ for all $\nu>\nu_c(g)$ and $\chi(g, \nu_c(g)) = \infty$, so our definition is consistent with the usual definition of $\nu_c(g)$. 
We also prove $\nu_c(g) \le 0$ in Proposition~\ref{prop:chi}. 
Limits of $\nu_c(g)$ are proved in Proposition~\ref{prop:limit-nu_c}. 

\begin{lemma} \label{lemma:nu_c}
Fix $g>0$. Consider $\nu\in \R$. 
\begin{enumerate}
\item[(i)] The map $\nu\mapsto \|Q(g, \nu)\|$ is continuous and strictly decreasing.
\item[(ii)] The map $\nu\mapsto \|Q(g, \nu)\|$ is convex. 
\item[(iii)] $\nu_c(g)$ exists and is unique. 
\end{enumerate}
\end{lemma}

\begin{proof}
(i) Since $\nu\in \R$, we have $k_0 > 0$, so \begin{equation}
\|Q\| = \sup_{\|f\|_2=1} \abs{ \langle Qf, f \rangle }
= \sup_{\|f\|_2=1, f\ge0} \langle Qf, f \rangle. \label{eqn:26-2}
\end{equation}
Consider $\nu, \hat \nu\in \R$. Denote $\hat  Q= Q(g, \hat \nu), \hat k_0(t,s) = \hat k_0(t,s;g, \hat \nu)$, and $ \hat p(t) = p(t; g, \hat \nu)$. For any $f\ge 0$, by the Cauchy--Schwarz inequality, \begin{align} \label{eqn:2.9}
\abs{ \langle Qf, f \rangle - \langle \hat Qf, f \rangle }
&= \left\lvert \int_0^\infty \int_0^\infty f(t)f(s) [ k_0(t,s) - \hat k_0(t,s) ] ds dt \right\lvert \\
&\le \|f\|_2^2  \big \| k_0(t,s) - \hat k_0(t,s)
\big \|_{L^2(dsdt)}.  \nonumber
\end{align}
Recall $k_0(t,s) = \sqrt{p(t)}\sqrt{p(s)} e^{-t}e^{-s} I_0(2\sqrt{st})$ has exponential decay by assumption~\eqref{A3}. 
Also notice $p(t) = e^{-g\phi(t) - \nu t}$ is monotone in $\nu$. 
Hence, as $\hat \nu\to \nu$, we have $\big \| k_0(t,s) - \hat k_0(t,s)
\big \|_{L^2(dsdt)}\to 0$ by the Dominated Convergence Theorem. 
Since we are taking supremum over $\|f \|_2=1$ to get $\|Q\|$, we obtain $\| \hat Q \| \to \| Q\|$. 

To prove that the map is strictly decreasing, consider $\nu > \hat \nu$. 
By Lemma~\ref{lemma:Q}(ii), $Q$ has an eigenvector $h$ with eigenvalue $\|Q \|$ that is continuous and strictly positive everywhere. We normalize it so that $\|h\|_2 = 1$. 
For each $t>0$, we have $p(t) = e^{-g\phi(t) - \nu t} < e^{-g\phi(t) - \hat\nu t} = \hat p(t)$, so $k_0(t,s) < \hat k_0(t,s)$ for all $t,s>0$. 
Since $h$ is positive and continuous, we get
\begin{align}
\| Q\| = \langle Qh, h \rangle 
&=\int_0^\infty \int_0^\infty h(t)h(s) k_0(t,s) ds dt \\
&< \int_0^\infty \int_0^\infty h(t)h(s) \hat k_0(t,s) ds dt
= \langle \hat Qh, h \rangle  \le \| \hat Q\|. \nonumber
\end{align}

\noindent
(ii) For each $f\ge0$, a direct calculation gives $\frac{\del^2}{\del\nu^2} \langle Qf, f \rangle \ge 0$,
so the map $\nu\mapsto \langle Qf, f \rangle$ is convex. 
By taking supremum \eqref{eqn:26-2}, the map $\nu\mapsto \|Q\|$ is convex too. 

\smallskip \noindent
(iii) As $\nu\to \infty$, $p(t) = e^{-g\phi(t) - \nu t} \to 0$ for all $t>0$. 
Similar to \eqref{eqn:2.9}, we have $\langle Qf, f \rangle\to 0$ as $\nu\to\infty$, uniformly in $\|f \|_2 =1$. It follows that $\| Q\| \to 0$ as $\nu\to \infty$. 
On the other hand, since $\nu\mapsto \|Q\|$ is convex and strictly decreasing, one of its subderivatives must be strictly negative. Since the corresponding subtangent line bounds the function from below, we have $\lim_{\nu\to -\infty} \|Q\| = \infty$. The critical $\nu_c(g)$ then exists by the Intermediate Value Theorem. Uniqueness follows from part (i). 
\end{proof}

\begin{lemma} \label{lemma:Q_analytic}
Let $g, \nu\in \C$ with $\Re(g)>0$. 
\begin{enumerate}
\item[(i)] 
For a fixed $g$, the map $\nu\mapsto Q(g, \nu)$ into the space of bounded linear operators on $L^2[0,\infty)$ is strongly analytic, with derivative given by \begin{equation}\label{def:dQdnu}
 \(\frac{\del Q}{\del\nu}\)f(t) = \int_0^\infty f(s) (-\half(t+s)) k_0(t,s) ds. 
\end{equation}
 
\item[(ii)] 
For a fixed $\nu$, the map $g\mapsto Q(g, \nu)$ is strongly analytic, with derivative given by \begin{equation}\label{def:dQdg}
\(\frac{\del Q}{\del g}\)f(t)  = \int_0^\infty f(s) (-\half(\phi(t)+\phi(s))) k_0(t,s) ds. 
\end{equation}
\end{enumerate}
\end{lemma}

\begin{proof}
(i) Fix some $g$ and $\nu$. Let $h\in\C$ with $0<\abs h \le1$. Define the operator \begin{align}
R_h = \frac{Q(g, \nu+h) - Q(g, \nu)}{h}. \end{align}
We will show $\lim_{h\to0} R_h = \frac{\del Q}{\del\nu}$ in the operator norm, where $\frac{\del Q}{\del\nu}$ (for now) is defined by the right-hand side of \eqref{def:dQdnu}. This will show $\nu\mapsto Q(g, \nu)$ to be strongly analytic at that point. 

Let $f_1, f_2\in L^2[0, \infty)$ be such that $\|f_1\|_2, \|f_2\|_2= 1$. 
Denote $p^*(t) = p(t \text{; } \Re(g), \Re(\nu))$ and $k_0^*(t,s) = k_0(t,s \text{; } \Re(g), \Re(\nu))$. 
By the Cauchy--Schwarz inequality, \begin{align}
&\left\lvert\left\langle  (R_h-\frac{\del Q}{\del\nu}) f_1, f_2\right\rangle \right\lvert   \\
&\le  
\int_0^\infty \int_0^\infty  \abs{ f_2(t) }\abs{ f_1(s) }
\left\lvert  \frac{e^{-\half h(t+s)}-1}{h}	- (-\half(t+s))\right\lvert 
k_0^*(t,s)dsdt \nonumber \\
& \le 
\|f_1\|_2\|f_2\|_2 \(\int_0^\infty \int_0^\infty \Bigg(
\left\lvert  \frac{e^{-\half h(t+s)}-1}{h}	- (-\half(t+s))\right\lvert 
k_0^*(t,s)\right)^2dsdt\Bigg)^{1/2}. \nonumber
\end{align}
Taking supremum over $f_1$ and $f_2$ gives 
\begin{align}
\left\| R_h-\frac{\del Q}{\del\nu}\right\| 
&= \sup_{\|f_1\|_2=\|f_2\|_2=1} \left\lvert\left\langle  (R_h-\frac{\del Q}{\del\nu}) f_1, f_2\right\rangle \right\lvert   \\
&\le \(\int_0^\infty \int_0^\infty \Bigg(
\left\lvert  \frac{e^{-\half h(t+s)}-1}{h}	- (-\half(t+s))\right\lvert 
k_0^*(t,s) \right)^2dsdt\Bigg)^{1/2}. \nonumber
\end{align}
Since $\abs h \le 1$, we have $\abs{ e^{ah} - 1 - ah }  
\le \sum_{n=2}^\infty \frac{ \abs{ ah }^n }{n!}  
\le \abs h^2 \sum_{n=2}^\infty \frac{ \abs{ a }^n }{n!}  
\le  \abs{h}^2 e^{\abs a}$ for all $a$.
We apply the inequality with $a = -\half (t+s)$, then \begin{align}
\left\lvert\frac{e^{-\half h(t+s)}-1}{h}	- (-\half(t+s))\right\lvert
\le \abs h e^{\half(t+s)}. 
\end{align}
By assumption~\eqref{A3}, $e^{\half(t+s)}k_0^*(t,s)$ still has exponential decay. The claim then follows from the Dominated Convergence Theorem. 

\smallskip \noindent
(ii) The same proof with $\abs h < \min \{ \Re(g), 1 \}$ works. This and assumption~\eqref{A3} guarantee enough decay to apply the Dominated Convergence Theorem. 
\end{proof}

\begin{remark}\label{remark:kato}
For $g>0$ and $ \nu\in \R$, the eigenvalue $\|Q(g, \nu)\|$ is simple and isolated by Lemma~\ref{lemma:Q}.
Since the operators $Q(g, \nu)$ are strongly analytic, the Kato--Rellich Theorem \cite[Theorem XII.8]{ReedSimon} implies that there exists a holomorphic function $\lambda(\nu;g)$ that agrees with $\norm{ Q(g, \nu) }$ when $\nu\in \R$. 
In particular, the derivative $\del_\nu \|Q(g, \nu)\|$ exists. 
Convexity and strict monotonicity (Lemma~\ref{lemma:nu_c}) then imply $ \del_\nu \|Q(g, \nu)\| <0$. 
Similarly, $\del_g \|Q(g, \nu)\|$ exists too. 

By assumptions~\eqref{A3}, there exists $t_0\ge0$ for which $\phi(t)>0$ for all $t> t_0$. 
Using this, the same proof of Lemma~\ref{lemma:nu_c} shows the map $g\mapsto \norm{ Q(g, \nu) }$ for a fixed $\nu \in \R$ to be convex and strictly decreasing. It follows that $\del_g \|Q(g, \nu)\|<0$. Then, by the Implicit Function Theorem, \begin{align} \label{eqn:2.21}
\frac {d\nu_c}{dg} = - \frac{ \del_g \|Q(g, \nu)\| }{ \del_\nu \|Q(g, \nu)\| }<0,
\end{align}
showing $\nu_c(g)$ to be strictly decreasing in $g$.
\end{remark}

\begin{proposition} \label{prop:limit-nu_c}
Assume $\nu_c(g) \le 0$ for all $g>0$ (which will be proved in Proposition~\ref{prop:chi}).
The critical parameter $\nu_c(g)$ is strictly decreasing in $g$, with limits \begin{align}
\lim_{g\to 0^+} \nu_c(g) = 0, \qquad 
\lim_{g\to \infty} \nu_c(g) = -\infty.
\end{align} 
\end{proposition}

\begin{proof}
The function $\nu_c(g)$ is strictly decreasing because $\nu_c'(g) < 0$ by \eqref{eqn:2.21}. It follows that limits of $\nu_c(g)$ exist. 
Since $\nu_c(g)\le 0$ by hypothesis, we have $\lim_{g\to 0^+} \nu_c(g) \le 0$. 

Suppose for the sake of contradiction that $\lim_{g\to 0^+} \nu_c(g) = \tilde \nu_0 < 0$. 
Let $\nu_0 = \max \{ \tilde \nu_0, -1 \}$. 
Since $\norm{ Q(g, \nu) } $ is decreasing in $\nu$ by Lemma~\ref{lemma:nu_c}(i), 
for any $g>0$ and $f\ne 0$ we have 
\begin{align} \label{eqn:2.23}
1 = \norm{ Q(g, \nu_c(g)) } \ge \norm { Q(g, \tilde \nu_0)  } \ge \norm { Q(g, \nu_0)  }
\ge \frac{ \norm{Q(g, \nu_0) f}_2 } {\norm {f}_2 }.
\end{align} 
We will pick $f$ so that the right-hand side $\to \infty$ as $g\to 0^+$.
Let $a>0$ be a parameter and let $f(t) = e^{-at}$. 
Since $\sqrt{p(t; g, \nu_0)}$ increases to $e^{-\half \nu_0 t}$ as $g\to 0^+$, by the Monotone Convergence Theorem and the Taylor series $I_0(z) = \sum_{m=0}^\infty \frac{1}{ m! m!} (\frac z 2)^{2m}$, 
\begin{align}
\lim_{g\to 0^+} (Q(g, \nu_0)f )(t)
&= \int_0^\infty f(s) e^{-\half \nu_0 t} e^{-\half \nu_0 s} e^{-t}e^{-s} I_0(2\sqrt{st})\ ds  \\
&= \sum_{m=0}^\infty \frac{ t^m} {m!} e^{-(\half \nu_0 + 1 ) t}   \int_0^\infty  \frac{ s^m} {m!} e^{-(a + \half \nu_0 + 1 ) s} ds  \nonumber \\
&=  \sum_{m=0}^\infty \frac{ t^m} {m!} e^{-(\half \nu_0 + 1 ) t}  \( a+ \half \nu_0 + 1 \)^{-m-1} \nonumber \\
&= \( a+ \half \nu_0 + 1 \)\inv \exp\Big\{ ( a+ \half \nu_0 + 1 )\inv t - (\half \nu_0 + 1 ) t \Big\}.  \nonumber
\end{align}
Observe that if 
\begin{align} \label{eqn:2.25}
( a+ \half \nu_0 + 1 )\inv  - (\half \nu_0 + 1 )  \ge 0, 
\end{align}
then $ \lim_{g\to 0^+} \norm{  Q(g, \nu_0)f  }_2 = \infty $ by monotone convergence. Since $\nu_0  \in (-2,0) $, equation \eqref{eqn:2.25} holds as a strict inequality when $a=0$. By continuity, \eqref{eqn:2.25} also holds for all small $a>0$. 
Picking one such $a$, then taking the $g\to 0^+$ limit of equation~\eqref{eqn:2.23} yields \begin{align}
1 \ge \frac{  \lim_{g\to 0^+} \norm{  Q(g, \nu_0)f  }_2 } { \norm f_2} = \infty,
\end{align} giving a contradiction. 

The $g\to \infty$ limit of $\nu_c(g)$ is easier. 
If $\nu_c(g) \to \nu_1 > -\infty$, 
since $\norm{ Q(g, \nu) } $ is decreasing in $\nu$, 
for any $g>0$ we have
\begin{align} 
1 = \norm{ Q(g, \nu_c(g)) } \le \norm { Q(g, \nu_1)  }. 
\end{align} 
As in the proof of Lemma~\ref{lemma:nu_c}(iii) (now with $g\to \infty$ instead if $\nu\to \infty$), the right-hand side converges to 0 as $g\to \infty$, giving a contradiction. 

This completes the proof.
\end{proof}

\subsection{Infinite-volume two-point function}
In this section, we prove $G_{ij}^N \to G_{ij}$ and give an inner product representation of $G_{ij}$. As a corollary, we obtain asymptotics of $G_{ij}$ as $\abs{ j-i } \to\infty$. 

\begin{proposition} \label{prop:Gij}
Let $g>0$. There exists $\eps = \eps(g) >0$ such that for $\Re(\nu) > \nu_c(g) - \eps$:
\begin{enumerate}
\item[(i)]
The limit $q = \lim_{N\to\infty} T^N[\sqrt p]$ exists in $L^2[0, \infty)$. Moreover, the convergence is locally uniform in $g$ and $\nu$, and $q(t; g, \nu)$ is continuous in $t, g$ and holomorphic in $\nu$;

\item[(ii)] 
For any integers $i \le j$, the two-point function $G_{ij}(g, \nu)$ is finite and is given by \begin{align} \label{eqn:2.28}
G_{ij}(g, \nu) = \lim_{N\to\infty} G_{ij}^N(g, \nu)
=  \langle Q^{j-i}(q) , \bar q  \rangle. 
\end{align}
\end{enumerate}
\end{proposition}

The proposition is proved via the contraction mapping principle. 
Recall $\nu_c(g) \in \R$ makes the operator norm $\|Q(g, \nu_c(g))\|=1$. 
The $Q^{j-i}$ in the representation \eqref{eqn:2.28} hints at different behaviors of $G_{ij}$ when $\norm Q<1$ and $\norm Q >1$ as $\abs{j-i} \to \infty$. 
\begin{corollary} \label{corollary:Gij}
Let $g>0$, $\nu\in \R$, and $\eps$ be given by Proposition~\ref{prop:Gij}. Then
\begin{enumerate}
\item[(i)]
If $\nu > \nu_c(g)$, $G_{ij}(g, \nu)$ decays exponentially as $\abs { j-i } \to \infty$. 

\item[(ii)]
If $\nu = \nu_c(g)$, $G_{ij}(g, \nu_c(g))$ converges to a non-trivial limit as $\abs { j-i } \to \infty$. 

\item[(iii)]
If $\nu_c(g) - \eps< \nu < \nu_c(g)$, $G_{ij}(g, \nu)$ diverges exponentially as $\abs { j-i } \to \infty$. 
\end{enumerate}
In particular, we have the critical exponent $\eta=1$ (see \eqref{exponent:eta}). 
\end{corollary}

\begin{proof}[Proof of Corollary~\ref{corollary:Gij}]
When $\nu\in \R$, by Lemma~\ref{lemma:Q}, $Q$ is compact, positive, self-adjoint, and $\norm Q$ is an isolated eigenvalue with a strictly positive eigenvector $h$.
Let $P_{\|Q\|}$ denote the projection into the eigenspace of $\norm Q$. 
If $P_{\|Q\|}(q) \ne 0$, then \eqref{eqn:2.28} and the spectral theorem implies \begin{align}
G_{ij}(g, \nu) \sim \|Q\|^{\abs{ j-i }} \| P_{\|Q\|}(q)\|_2^2, \qquad \abs{ j-i }\to \infty, \end{align}
which gives the desired result. Hence, it is sufficient to prove $P_{\|Q\|}(q) \ne 0$. 

Recall the definition of the operator $T$ in \eqref{def:T}. 
Since $Tf(t) \ge \sqrt{p(t)}e^{-t} \ge 0$ when $f \ge 0$, 
a simple induction yields that 
$T^N[\sqrt p](t) \ge \sqrt{p(t)}e^{-t}$ for all $N$. 
It then follows from the $L^2$ convergence in Proposition~\ref{prop:Gij}(i) that
\begin{align} 
q(t) \ge \sqrt{p(t)}e^{-t} 
\end{align}
for almost every $t$, 
which implies $\langle q, h \rangle > 0$ and $P_{\|Q\|}(q) \ne 0$. 
This concludes the proof. 
\end{proof}

\subsubsection{Proof of Proposition~\ref{prop:Gij}(i)}
\emph{Step 1.}
We first consider a fixed pair of parameters $g>0$, $\nu\in\R$,
and define the auxiliary linear operator
\begin{align}
Af(t) = \int_0^\infty f(s) k_1(t,s)\sqrt{\tfrac{t}{s}} ds, 
\end{align}
so that $Tf(t) = \sqrt{p(t)}e^{-t} + Af(t)$. 
We will view $A$ as an operator on a weighted $L^r$ space and show that it is a contraction when $\nu\ge \nu_c(g)$.

For $1 < r < 2$, we define $w(t) = t^{-r/2}$ and the measure $d\mu = w(t) dt$ on $(0, \infty)$. Let $r' \in (2, \infty)$ denote the H\"older conjugate of $r$. 
Since $\sqrt{ p(t) }$ is bounded near $t=0$ and has exponential decay, 
$\sqrt{p(t)}$ belongs to  $L^r(\mu) = L^r_\mu$ for all $r < 2$.
For a function $h\in L^{r'}_\mu$, by definition of $r'$, we have
\begin{align} \label{eq:change_of_measure}
\int_0^\infty \abs{ h(t) }^{r'} d\mu(t)
= \int_0^\infty \abs{ h(t) }^{r'} \frac{1}{t^{(r-1)r'/2}} dt
= \int_0^\infty \( \abs{ h(t) }  t^{ (1-r)/2 } \)^{r'} dt,
\end{align}
so $h\in L^{r'}_\mu$ if and only if $ h(t) t^{ (1-r)/2 } \in L^{r'}_\dt$, 
with $\norm h_{L^{r'}_\mu} = \norm{ h(t) t^{ (1-r)/2 } }_{L^{r'}_\dt}$. 
Similarly, for $f\in L^r_\mu$, we have $f \in L^r_\mu$ if and only if $ f(t) t^{-1/2} \in L^r_\dt$, with $\norm f _{L^r_\mu} = \norm{ f(t) t^{-1/2} }_{ L^r_\dt}$. 
We view $A$ as an operator from $L^r_\mu \to L^r_\mu$, 
and define another operator $B:L^r_\dt \to L^r_\dt$ by
\begin{align}
Bf(t) = \int_0^\infty f(s) k_1(t,s)ds. 
\end{align}
Both $A$ and $B$ are bounded operators because $k_1(\cdot, \cdot)$ has exponential decay by assumption~\eqref{A3}. 
For $f\in L^r_\mu$ and $h\in L^{r'}_\mu$, observe that
\begin{align}
\langle Af, h \rangle_\mu
&= \int_0^\infty  \overline{h(t)} \int_0^\infty f(s) k_1(t,s)\sqrt{\tfrac{t}{s}}\ ds d\mu(t)  \\
&= \int_0^\infty \frac{ \overline{h(t)}\sqrt{t}}{t^{r/2}} \int_0^\infty \frac{f(s)}{\sqrt{s}} k_1(t,s)\ ds dt \nonumber \\
&= \left\langle B\( f(t) t^{-1/2}   \), h(t) t^{(1-r)/2}  \right\rangle_{dt}, \nonumber
\end{align}
so \begin{align}
\|A\|_{L^r_\mu}
&= \sup \{ \abs{ \langle Af, h \rangle_\mu } : 
	\| f \|_{L^r_\mu} = \|h\|_{L^{r'}_\mu} =1 \} \\
&= \sup \{ \abs{ \langle B\hat f, \hat h \rangle_\dt } : 
	\| \hat f \|_{L^r_\dt} = \|\hat h\|_{L^{r'}_\dt} =1 \}
= \|B\|_{L^r_\dt}.\nonumber
\end{align}
Thus, to prove that $A$ is a contraction, it suffices to prove $ \|B\|_{L^r_\dt} < 1$ for some $1<r<2$. For simplicity, we write $\|B\|_r = \|B\|_{L^r_\dt}$. 

For $r=2$, the methods of Section~\ref{section:Qproperty} show that $B:L^2_\dt \to L^2_\dt$ is compact, self-adjoint, and has a continuous, strictly positive eigenvector $h$ with eigenvalue $\|B\|_2$, which we normalize to have $\|h\|_2=1$. Since the modified Bessel functions satisfy $I_1<I_0$ pointwise, 
\begin{align}
\|B\|_2 = \langle Bh, h \rangle < \langle Qh, h \rangle \le \|Q\|.
\end{align}
Then by the Riesz--Thorin interpolation theorem, for $1 < r \le 2$, we have
\begin{align}
\| B\|_r  \le \| B\|_1^{1-\theta(r)} \| B\|_2^{\theta(r)} 
< \| B\|_1^{1-\theta(r)} \| Q \|^{\theta(r)} 
\end{align}
for some continuous function $\theta(r) \in (0,1]$ with $\theta(2) = 1$. 
Since $\norm B_1 < \infty$, continuity implies that
$\|B\|_r < \| Q\|$ when $r$ is sufficiently close to 2. 
Hence, for $\nu\ge \nu_c(g)$ (which satisfies $\|Q(g, \nu)\|\le 1$) and $r$ slightly less than 2, we get $\|A\|_{L^r_\mu} = \| B\|_r < 1$, showing $A$ to be a contraction in $L^r_\mu$.

\smallskip \noindent
\emph{Step 2.} We now allow $\nu$ to be complex. 
Let $g_0>0,\ \nu_0\in \R$ be such that $\| Q(g_0, \nu_0) \| \le 1$. 
By Step 1, there exists some $r_0 \in (1,2)$ for which $\| B(g_0, \nu_0) \|_{r_0} <1$. 
Since the map $(g, \nu)\mapsto \| B(g, \nu) \|_{r_0}$ is continuous, there exists $\delta = \delta(g_0, \nu_0) >0$ such that $\| B(g, \nu) \|_{r_0}<1$ for all $\abs { g-g_0 }\le\delta$, $\abs { \nu-\nu_0 } \le \delta$. 
We first show uniform $L^{r_0}_\mu$ convergence for $(g,\nu)$ in the tube
\begin{align}
S_0= \{ (g, \nu) \in (0, \infty) \times \C : \abs{g - g_0} \le \delta,\ \abs{\Re(\nu) - \nu_0} \le \delta \}, 
\end{align}
where $d\mu(t) = t^{- r_0/2} dt$,
and then we prove $L^2_\dt$ convergence using H\"older's inequality. 
The $\eps$ in the statement of the proposition is given by $\eps(g_0) = \delta(g_0, \nu_c(g_0))$.

We extend the definitions of $T$ and $A$ to include dependence on $g$ and $\nu$. 
Let 
\begin{align}
\mathcal X_0 = \{  \psi : (0, \infty)\times S_0 \to \C 
\mid \norm{ \psi }_{\mathcal X_0} = 
\sup_{(g, \nu)\in S_0} \| \psi (\cdot, g, \nu) \|_{L^{r_0}_\mu }<\infty\} .
\end{align} 
We define operator $\hat T:\mathcal X_0\to \mathcal X_0$ and linear operator $\hat A:\mathcal X_0\to \mathcal X_0$ by 
\begin{align}
(\hat A\psi)(t, g, \nu) &= \int_0^\infty \psi(s, g, \nu) k_1(t,s;g,\nu)\sqrt{\tfrac{t}{s}} ds, \\
(\hat T\psi)(t, g, \nu) &= \sqrt{p(t; g, \nu)} e^{-t} + (\hat A\psi)(t, g, \nu). 
\end{align}
These are indeed well-defined operators into $\mathcal X_0$, because all $(g, \nu)\in S_0$ are dominated by the left real boundary point $g^* = g_0-\delta$, $\nu^* = \nu_0-\delta$, in the sense that 
\begin{align} \label{eq:p*}
\abs{\sqrt{p(t; g, \nu)} } = \abs{e^{- \half g\phi(t)- \half \nu t} } 
\le e^{ - \half g^*\phi(t) - \half \nu^* t } = \sqrt{p^*(t)},
\end{align} 
where the last equality defines $p^*(t)$. 
Writing $A^* = A(g^*, \nu^*)$, we have
\begin{align}
\abs{(\hat A\psi)(t, g, \nu)} \le \big(A^*\abs{\psi(\cdot, g, \nu)}\big)(t)
	\quad \forall t, 
\end{align}
so \begin{align}
\norm{ \hat A\psi (\cdot, g, \nu)  }_{L^{r_0}_\mu}
&\le \big\| A^*\abs{\psi(\cdot, g, \nu)} \big\|_{L^{r_0}_\mu} 
\le \| A^* \|_{L^{r_0}_\mu} \| \psi(\cdot, g, \nu) \|_{L^{r_0}_\mu}
	~~\forall\, (g, \nu)\in S_0. 
\end{align}
Taking supremum yields 
\begin{align}
\sup_{(g, \nu)\in S_0}\| \hat A\psi (\cdot, g, \nu) \|_{L^{r_0}_\mu}
\le \| A^* \|_{L^{r_0}_\mu} \sup_{(g, \nu)\in S_0} \| \psi (\cdot, g, \nu) \|_{L^{r_0}_\mu}, 
\end{align}
which gives the operator norm $\| \hat A \| \le \| A^* \|_{L^{r_0}_\mu} = \| B(g^*, \nu^*) \|_{r_0}<1$, showing $\hat A$ to be a contraction on $\mathcal X_0$. 
Since \begin{align} \label{eq:T_diff}
\hat T\psi_1 - \hat T \psi_2 
= \hat A\psi_1 - \hat A \psi_2 
\quad \forall \psi_1, \psi_2 \in \mathcal X_0,
\end{align} 
$\hat T$ is also a contraction. It follows that $q = \lim_{N\to\infty} \hat T^N[\sqrt p]$ exists in $\mathcal X_0$ and is a fixed point of $\hat T$. 

We next show that the convergence also happens in $L^2_\dt$, uniformly in $(g, \nu)\in S_0$. 
Let $k_1^*(t,s) = k_1(t,s; g^*, \nu^*)$. 
Since $\hat Tq = q$, by \eqref{eq:T_diff}, \eqref{eq:p*}, H\"older's inequality, and \eqref{eq:change_of_measure}, 
\begin{align}
&\left\lvert \( q -  \hat T^{N+1}[\sqrt p] \)(t, g, \nu)\right\lvert  
= \left\lvert \hat A\( q -  \hat T^{N}[\sqrt p] \)(t, g, \nu)\right\lvert  \\
&\qquad\qquad \le \int_0^\infty \left\lvert\(  q -  \hat T^{N}[\sqrt p]   \)(s, g, \nu)\right\lvert
 k_1^*(t,s) \sqrt{\tfrac{t}{s}}s^{r_0/2} d\mu(s) \nonumber \\
&\qquad\qquad \le \left\| \(q -  \hat T^{N}[\sqrt p] \)(\cdot, g, \nu) \right\|_{L^{r_0}_\mu} 
\Big\|     k_1^*(t,s) \sqrt{\tfrac{t}{s}}s^{r_0/2}       \Big\|_{L^{r_0'}_{d\mu(s)}} \nonumber\\
&\qquad\qquad = 
\left\| \(q -  \hat T^{N}[\sqrt p] \)(\cdot, g, \nu) \right\|_{L^{r_0}_\mu} 
\big\|     k_1^*(t,s) \sqrt{t}  \big\|_{L^{r_0'}_{ds}} . \nonumber
\end{align}
Taking supremum over $(g, \nu)\in S_0$, and then the $L^2_\dt$ norm, we get
\begin{align} \label{eq:L2_conv} 
\bigg\lVert
\sup_{g, \nu} \left\lvert \( q -  \hat T^{N+1}[\sqrt p] \)(t, g, \nu)\right\lvert  
\bigg\rVert _{ L^2_\dt }
\le   \left\| q -  \hat T^{N}[\sqrt p] \right\|_{\mathcal X_0} 
	\left\lVert \big\|     k_1^*(t,s) \sqrt{t}  \big\|_{L^{r_0'}_{ds}} \right\rVert_{L^2_\dt}.
\end{align}
The first term goes to 0 by the convergence in $\mathcal X_0$, and the second term is finite because $k_1^*$ has exponential decay by assumption~\eqref{A3}. 
This proves that $T^N[\sqrt p] \to q$ in $L^2_\dt$.
Also, replacing the $L^2_\dt$ norm in \eqref{eq:L2_conv} by a supremum over $t\ge 0$ shows that $T^N[\sqrt p] \to q$ uniformly in $t\ge 0$ and in $(g, \nu) \in S_0$.
It follows that $q(t; g, \nu)$ is continuous and is holomorphic in $\nu$
(as all $T^N[\sqrt p]$ are holomorphic in $\nu$). 
This concludes the proof of Proposition~\ref{prop:Gij}(i). 
\qed

\subsubsection{Proof of Proposition~\ref{prop:Gij}(ii)}
Let $i\le j$.
The second equality in the claim \eqref{eqn:2.28} follows directly from Proposition~\ref{prop:Gij^N}, part (i), and the continuity of $Q$. 
To prove the first equality, we first prove $P_{ij}^N(g, T) \to P_{ij}(g, T)$. 
This argument is adapted from \cite{BBS2015}. Notice if a walk never reaches $\pm N$, then it contributes the same to $E_i^N[\cdot]$ and $E_i[\cdot]$.  Therefore, a walk starting at $i$ must make at least $\min\{\abs{N-i}, \abs{-N-i}\}$ steps to make a difference. Using $L_{T,x}\ge0$ for all $x$, 
\begin{align} 
\abs{ P_{ij}^N(g, T) - P_{ij}(g,T)}
&\le 2 P_i( X \text{ hits one of }\pm N \text{ by time }T) \\
&= 2P(M_T \ge \min\{\abs{N-i}, \abs{-N-i}\}),  \nonumber
\end{align}
where $M_T$ is a Poisson process with rate $2$. 
For a fixed $T$, this probability converges to 0 as $N\to\infty$. 

Since $P_{ij}^N\le 1$ by definition, when $\nu \in \R$  and $\nu > 0$, 
dominate convergence implies that
\begin{align} \label{eqn:51}
G_{ij}(g, \nu) 
= \lim_{N\to\infty} \int_0^\infty P_{ij}^N(g, T) e^{-\nu T}dT
= \lim_{N\to\infty} G_{ij}^N(g, \nu)
= \langle Q^{j-i}(q) , \bar q  \rangle .
\end{align}
Note both sides of \eqref{eqn:51} are finite and holomorphic in $\nu$ for $\Re(\nu)> \nu_c(g) - \eps$, as
\begin{align} 
\abs{ G_{ij}(g, \nu) } 
\le G_{ij}(g, \Re(\nu) )
\le \lim_{N\to\infty} G_{ij}^N(g, \Re(\nu) )
=  \langle Q^{j-i}(q) , \bar q  \rangle \Big\rvert_{g, \Re(\nu)} 
\end{align}
by Fatou's lemma. 
By the uniqueness of analytic continuation, \eqref{eqn:51} must hold for all $\Re(\nu)> \nu_c(g) - \eps$, giving the first equality of \eqref{eqn:2.28}.
\qed

\subsection{Susceptibility and correlation length} 
\label{section:susceptibility}
In this section, we prove results about the susceptibility and correlation lengths. 

\begin{proposition} \label{prop:chi}
Let $g > 0$, $\nu_c = \nu_c(g)$, $\nu \in \C$ with $\Re(\nu) > \nu_c$,
and $q$ be given by Proposition~\ref{prop:Gij}. Then: 
\begin{enumerate}
\item[(i)]
The one-sided susceptibility $\chi_+(g, \nu)$ is given by 
\begin{align} \label{eqn:chi}
\chi_+(g, \nu) = \sum_{j=1}^\infty G_{0j}(g, \nu) 
= \left \langle Q(1-Q)\inv(q) , \bar q \right \rangle,
\end{align}
and there is a constant $\bar u = \bar u(g)>0$ such that 
\begin{align}  \label{eqn:chi_asym}
\chi_+(g, \nu)
\sim  \frac{\bar u}{\nu-\nu_c}
\(-\frac{\del \|Q\|}{\del\nu}\Big\lvert_{\nu=\nu_c(g)}\)\inv
, \qquad \nu\to \nu_c^+. 
\end{align}

\item[(ii)]
$\nu_c(g) \le 0$. 
\end{enumerate}
\end{proposition}
We have isolated a constant $\bar u$ in the residue of $\chi_+$ in~\eqref{eqn:chi_asym}, because $\bar u$ shows up in the calculation of correlation lengths also. 
By symmetry, the two-sided susceptibility is given by $\chi(g, \nu) = 2\chi_+(g, \nu) + G_{00}(g, \nu)$. Since $G_{00}$ is regular as $\nu\to \nu_c$ by Proposition~\ref{prop:Gij}(ii), we obtain Theorem~\ref{theorem:critical_exponent}(ii). 

\begin{proof}[Proof of Proposition~\ref{prop:chi}]
(i) Let $\Re(\nu) > \nu_c$, so $\|Q(g, \nu) \| \le \| Q(g, \Re(\nu)) \| < 1$. 
Using Proposition~\ref{prop:Gij}(ii), we have
\begin{align}
\abs{ \chi_+(g, \nu) } \le \sum_{j=1}^\infty \abs { \langle Q^{j}(q) , \bar q  \rangle }
\le \sum_{j=1}^\infty \|Q\|^j \| q\|_2^2 < \infty.
\end{align} 
Since the convergence is locally uniform in $\nu$, $\chi_+$ is holomorphic in $\nu$. 
When such $\nu$ is real, by the Monotone Convergence Theorem and holomorphic functional calculus, we have
\begin{align}
\chi_+(g, \nu) = \sum_{j=1}^\infty \langle Q^{j}(q) , \bar q  \rangle
=  \left \langle \sum_{j=1}^\infty Q^{j}(q) , \bar q  \right\rangle
= \left \langle Q(1-Q)\inv(q) , \bar q \right \rangle. 
\end{align} 
It then follows from the uniqueness of analytic continuation that the equality holds for all $\Re(\nu)>\nu_c$.

To determine the divergence of $\chi_+$ as $\nu\to \nu_c^+$, we restrict to $\nu\in\R$ so that $Q$ is self-adjoint. By the spectral theorem, we can decompose $q= v + w$ where $v = P_{\|Q\|}(q)$ is the projection of $q$ into the eigenspace $E_{\|Q\|}$, and $w\in (E_{\|Q\|})^\perp$. 
Since $\| Q\|$ is a simple eigenvalue, by the Kato--Rellich Theorem \cite[Theorem XII.8]{ReedSimon}, this decomposition is analytic around $\nu=\nu_c$ (where $\|Q\|=1$), and there exists $\delta>0$ such that the spectrum of $Q(g, \nu)$ intersects $\{ \lambda\in \C \mid \abs{\lambda-1} < \delta\}$ at exactly one point, for all $\nu$ sufficiently close to $\nu_c$. In particular, this implies that all other eigenvalues of $Q$ are bounded by $1-\delta$ as $\nu\to \nu_c$. Hence, using that $q$ is continuous at $\nu_c$, we have
\begin{align}
\chi_+(g, \nu) 
&= \left \langle Q(1-Q)\inv(q) , \bar q \right \rangle  \\
&= \left \langle Q(1-Q)\inv(v) , \bar v \right \rangle + \left \langle Q(1-Q)\inv(w) , \bar w \right \rangle\nonumber \\
&= \frac{\|Q\| }{1-\| Q\|} \left \| P_{\|Q\|}(q) \right \|_2^2 + \left \langle Q (1-Q)\inv(w) , \bar w \right \rangle\nonumber \\
&\sim \frac{1}{1-\| Q\|} \left \| P_{1}(q\lvert_{\nu_c}) \right \|_2^2 + O(1)
, \qquad\qquad\qquad \nu\to \nu_c^+.\nonumber
\end{align}
Since $\frac{\del \|Q\|}{\del\nu}\big\lvert_{\nu=\nu_c} = \lim_{\nu \to \nu_c } \frac{ \norm Q - 1}{\nu-\nu_c}$,
we get \begin{align}
\chi_+(g, \nu) \sim
\frac{1}{\nu-\nu_c}
\(-\frac{\del \|Q\|}{\del\nu}\Big\lvert_{\nu=\nu_c}\)\inv
\left \| P_{1}(q\lvert_{\nu_c}) \right \|_2^2 
, \qquad \nu\to \nu_c^+, 
\end{align} 
which is the desired result with $\bar u = \left \| P_{1}(q\lvert_{\nu_c}) \right \|_2^2 $. 

\smallskip \noindent
(ii) Since $\sum_{j=1}^\infty P_{0j} \le 1$, when $\nu > 0$, we have
\begin{align}
\chi_+(g, \nu) = \sum_{j=1}^\infty G_{0j}(g, \nu) 
\le \int_0^\infty 1\cdot e^{-\nu T}dT < \infty
\end{align}
by the Monotone Convergence Theorem. 
But $\chi_+(g, \nu_c(g)) = \infty$ by part (i), so we must have $\nu_c(g)\le 0$. 
\end{proof}

The same method applies to correlation lengths of integer orders. 
\begin{corollary} \label{corollary:correlation}
Let $g > 0$, $\nu_c = \nu_c(g)$, $\nu \in \C$ with $\Re(\nu) > \nu_c$,
and $q$ be given by Proposition~\ref{prop:Gij}. Then 
\begin{align} \label{eqn:2.54}
\sum_{j=1}^\infty jG_{0j}(g, \nu) 
= \left \langle Q(1-Q)^{-2}(q) , \bar q \right \rangle ,
\end{align}
and \begin{align}
\sum_{j=1}^\infty j G_{0j}(g, \nu) 
\sim \frac{\bar u}{(\nu-\nu_c)^2}
\(-\frac{\del \|Q\|}{\del\nu}\Big\lvert_{\nu=\nu_c(g)}\)^{-2}
, \qquad \nu\to \nu_c^+,
\end{align}
with the same constant $\bar u$ as in Proposition~\ref{prop:chi}. 

In general, for any $k\in \N$, we have
\begin{align} \label{eqn:j^k_asymp}
\sum_{j=1}^\infty j^k G_{0j}(g, \nu) 
\sim \frac{\bar u \cdot k!}{(\nu-\nu_c)^{k+1}}
\(-\frac{\del \|Q\|}{\del\nu}\Big\lvert_{\nu=\nu_c(g)}\)^{-(k+1)}
, \qquad \nu\to \nu_c^+. 
\end{align}
Thus, using symmetry, the correlation length of order $k$ (defined in \eqref{exponent:nu_k}) satisfies \begin{align}
\xi_k(g, \nu) 
\sim \frac{(k!)^{1/k}}{\nu-\nu_c}
\(-\frac{\del \|Q\|}{\del\nu}\Big\lvert_{\nu=\nu_c(g)}\)\inv
, \qquad \nu\to \nu_c^+, \end{align}
which gives the critical exponent $\nu_k=1$ . 
\end{corollary}

\begin{proof}
For $\Re(\nu) > \nu_c$, by the same argument as for Proposition~\ref{prop:chi}, \begin{align}
\sum_{j=1}^\infty j G_{0j}(g, \nu) 
= \sum_{j=1}^\infty j \langle Q^{j}(q) , \bar q  \rangle
=  \left \langle \sum_{j=1}^\infty j Q^{j}(q) , \bar q  \right\rangle
= \left \langle Q(1-Q)^{-2}(q) , \bar q \right \rangle. 
\end{align} 
The proof of the asymptotic formula is analogous. 
For $k\ge 2$, we calculate $\sum_{j=1}^\infty j^k Q^j$ by differentiating the geometric series, then we use the same argument. 
\end{proof}

\section{Time asymptotics} \label{section:laplace}

In this section, we prove asymptotic results as $T \to \infty$. We first prove a general Tauberian theorem, which utilizes analyticity properties of the Laplace transform on the boundary of the region of convergence. Then we apply the Tauberian theorem to prove Therorem~\ref{theorem:moments}.

For a function $f:[0, \infty) \to \R$, we define its Laplace transform $\Lcal f$ to be the complex-valued function \begin{align}
\Lcal f(z) = \int_0^\infty f(T) e^{-zT}dT. \end{align}

\subsection{Tauberian theorem}
The following is our Tauberian theorem. 

\begin{theorem} \label{theorem:tauber0}
Let $k\in \N_0$. Suppose $f:[0, \infty) \to[0, \infty)$ is differentiable and let its derivative be decomposed as $f' = \alpha_+ - \alpha_-$ with $ \alpha_\pm \ge 0$. 
Suppose each of the Laplace transforms $\Lcal f(z), \Lcal \alpha_\pm(z)$ converges for $\Re(z)>0$ and can be extended to a meromorphic function on an open set containing the closed half-plane $\{ z\in \C \mid \Re(z)\ge 0\}$. 
Suppose $\Lcal f(z)$ has a unique pole of order $k+1$ at $z=0$, and each of $\Lcal \alpha_\pm(z)$ either has a unique pole of order $\le k+1$ at $z=0$ or is holomorphic. 
If $\lim_{z\to 0} z^{k+1}\Lcal f(z) = C >0$, then \begin{align}
\lim_{T\to \infty} \frac{f(T)}{T^k} = \frac{C}{k!}. \end{align}
\end{theorem}

The main tools to prove our Tauberian theorem are the following two theorems. The first is from \cite[Theorem III.9.2]{Korevaar2004} and the second is from \cite[Theorem 4.1]{Korevaar1954}. 
\begin{theorem} \label{theorem:tauber1}
Let $\alpha(t)=0$ for $t<0$, be bounded from below for $t\ge 0$ and be such that the Laplace transform $G(z) = \Lcal\alpha(z)$ exists for $\Re(z)>0$. 
Suppose that $G(\cdot)$ can be analytically continued to an open set containing the closed half-plane $\{ z\in \C \mid \Re(z)\ge0\}$. Then the improper integral $\lim_{T\to\infty}\int_0^{T} \alpha(t)dt$ exists and equals $G(0)$. 
\end{theorem}
\begin{theorem} \label{theorem:tauber2}
Let $a(t)$ be integrable over every finite interval $(0, T)$, and let $\Lcal a(z)$ be convergent for $z>0$. Suppose $\Lcal a(z)$ can be extended analytically to a neighborhood of $z=0$. Finally, suppose that \begin{align}
a(t) \ge -\psi(t) \qquad (t>0),\end{align}
where $\psi(t)$ is continuous for $t>0$ and of the form $t^\gamma L(t)$, $L(t)$ slowly oscillating\footnote{$ L:(0, \infty) \to (0, \infty)$ is said to be slowly oscillating if it is continuous and $L(ct)/L(t)\to 1$ as $t\to\infty$ for every $c>0$. }. Then \begin{align}
\left\lvert\int_0^T a(t)dt - \lim_{z\to0} \Lcal a(z) \right\lvert= O(\psi(T)), \qquad T\to \infty. \end{align}
\end{theorem}

Notice Theorem~\ref{theorem:tauber2} does not assume $\psi(T) \to 0$, so the conclusion is different from that of Theorem~\ref{theorem:tauber1}. 
Under the hypotheses of Theorem~\ref{theorem:tauber1}, we can take the $\psi(T)$ in Theorem~\ref{theorem:tauber2} to be a constant function, then Theorem~\ref{theorem:tauber2} only gives that $ \int_0^T \alpha - G(0) = O(1)$, which is weaker than the conclusion of Theorem~\ref{theorem:tauber1}. Nevertheless, the flexibility of Theorem~\ref{theorem:tauber2} is that we can take $\psi(T)$ to be, for example, polynomials, and the consequent polynomial upper bound is sufficient for our purposes.

\begin{proof}[Proof of Theorem~\ref{theorem:tauber0}]
The framework of the proof is as follows. We will first use Theorem~\ref{theorem:tauber2} on a modification of $\alpha_\pm$ to prove $f(T) = O(T^k)$. Then we will use Theorem~\ref{theorem:tauber1} on a different modification of $f$ and $\alpha_\pm$ to show that $\lim_{T\to\infty} \frac{f(T)}{T^k} $ exists. Finally, we use the Hardy--Littlewood Tauberian theorem to identify the limit. 

By the assumptions on $\Lcal \alpha_\pm$, there are $A_j, B_j\in \R$ such that \begin{align}
\Lcal \alpha_+ (z) = \sum_{j=1}^{k+1} \frac{A_j}{z^j} + O(1), \qquad
\Lcal \alpha_- (z) = \sum_{j=1}^{k+1} \frac{B_j}{z^j} + O(1)
\end{align} as $z\to 0$. We claim $A_{k+1} = B_{k+1}$. 
This is because for $z>0$, integration by parts gives\footnote{For the boundary term, the existence of the limit $\lim_{T\to\infty} f(T)e^{-zT}$ follows from the existence of the Laplace transforms $\Lcal f(z)$ and $\Lcal f'(z)$.} $\Lcal[f'](z) = z\Lcal[f](z) - f(0)$. 
By assumption, $\Lcal f$ has a pole of order $k+1$ at 0, so $\Lcal f'(z)$ has a pole of order $k$ at 0. The relation $f' = \alpha_+ - \alpha_-$ then forces $A_{k+1} = B_{k+1}$. 

We subtract polynomials from $\alpha_\pm$ so that the Laplace transforms of the resultant functions no longer have poles. Let
\begin{align} \label{def:alpha-tilde}
\tilde \alpha_+(T) = \alpha_+(T) - \sum_{j=1}^{k+1} A_j \frac{T^{j-1}}{(j-1)!}, \qquad
\tilde \alpha_-(T) = \alpha_-(T) - \sum_{j=1}^{k+1} B_j \frac{T^{j-1}}{(j-1)!} .
\end{align}
Since $\Lcal[T^{j-1}/(j-1)!](z) = z^{-j}$, $\Lcal \tilde \alpha_\pm(z)$ are regular as $z\to 0$. From the assumptions on $\Lcal \alpha_\pm$, we get that $\Lcal \tilde \alpha_\pm$ extend analytically to an open set containing the full closed half-plane $\{ z\in\C \mid \Re(z) \ge 0\}$. 

We focus on $\alpha_+$, the argument for $\alpha_-$ is analogous. 
Since $\alpha_+(T) \ge 0$, we have $\tilde \alpha_+(T) \ge - \sum_{j=1}^{k+1} A_j \frac{T^{j-1}}{(j-1)!}$. 
To apply Theorem~\ref{theorem:tauber2}, we define \begin{align}
\psi(T) = T^k L(T)
= T^k \max\Big\{ \frac{ \abs{A_k} }{ k!} +1, \frac{1}{T^k} \Big\lvert \sum_{j=1}^{k+1} A_j \frac{T^{j-1}}{(j-1)!} \Big\rvert \Big\} , \end{align}
so $\alpha_+(T) \ge -\psi(T)$. 
For $T$ large enough, we have $L(T) = \frac{ \abs { A_k } }{ k!} +1$, so $L$ is slowly oscillating.  
Therefore, by Theorem~\ref{theorem:tauber2}, 
\begin{align}
\left\lvert \int_0^T \tilde \alpha_+(t)dt + \text{const} \right\lvert = O(T^k), \qquad T\to \infty.
\end{align} 
The same equation holds for $\tilde \alpha_-$ in the place of $\tilde \alpha_+$. We subtract the two equations and use the definition of $\tilde \alpha_\pm$, $f' = \alpha_+ - \alpha_-$, and $A_{k+1} = B_{k+1}$, to get 
\begin{align}
\bigg\lvert \int_0^T  \Big( f'(t) - \sum_{j=1}^k (A_j - B_j) \frac{t^{j-1}}{(j-1)!} \Big)dt + \text{const} \bigg\lvert = O(T^k), \qquad T\to \infty.
\end{align} 
Since the polynomial term integrates to $O(T^k)$, we get 
$f(T) = f(0) + \int_0^T f' =  O(T^k)$ as $T\to \infty$. 

Next, we will prove that $\lim_{T\to\infty} \frac{f(T)}{T^k} = \lim_{T\to\infty} \frac{f(T)}{(T+1)^k}$ exists. By the Fundamental Theorem of Calculus, \begin{align} \label{eqn:3.10}
\frac{f(T)}{(T+1)^k}
= f(0) + \int_0^T  \frac{ f'(t) }{ (t+1)^k} dt  -  k  \int_0^T \frac{ f(t) }{ (t+1)^{k+1}} dt. 
\end{align}
We first calculate the limit of the second integral using Theorem~\ref{theorem:tauber1}. 
Since $\Lcal f$ has a pole of order $k+1$, there are $C_j\in \R$ such that \begin{align}
\Lcal f(z) = \sum_{j=1}^{k+1} \frac{C_j}{z^j} + O(1), \qquad z\to 0,  \end{align} where $C_{k+1} = C$ in the hypotheses. 
Analogous to $\tilde \alpha_\pm$, we define \begin{align} 
\tilde f(T) = f(T) - \sum_{j=1}^{k+1} C_j \frac{T^{j-1}}{(j-1)!}, \end{align}
then $\Lcal \tilde f$ extends analytically to an open set containing the closed half-plane, but $\tilde f$ is no longer bounded from below. 
To fix this, we apply Theorem~\ref{theorem:tauber1} to $\tilde f(T) / (T+1)^{k+1}$, which is bounded from below. 
For $z>0$, the Laplace transform $\Lcal [ f(T) / (T+1)^{k+1} ](z)$ exists by domination, and $\Lcal [ \tilde f(T) / (T+1)^{k+1} ](z)$ exists by linearity. 
By a simple induction on $k$, 
\begin{align} \label{eqn:3.13}
\Lcal \Big[ \frac{\tilde f(T)}{(T+1)^{k+1}} \Big](z)
= e^z \int_z^\infty ds_1      \int_{s_1}^\infty ds_2\, \dots
\int_{s_{k}}^\infty ds_{k+1}  e^{-s_{k+1}} \Lcal[\tilde f ](s_{k+1}). \end{align}
This equation also holds for complex $z$ with $\Re(z)>0$ because $\Lcal \tilde f$ is analytic and the open half-plane is simply connected. 
Also, since $\Lcal \tilde f$ can be extended analytically to the closed half-plane, \eqref{eqn:3.13} extends $\Lcal [ \tilde f(T) / (T+1)^{k+1} ]$ analytically to the closed half-plane as well. 
Thus,  by Theorem~\ref{theorem:tauber1}, \begin{align}
\lim_{T\to \infty} \int_0^T \frac{\tilde f(t)}{(t+1)^{k+1}} dt
= \lim_{z\to0} \Lcal \Big[ \frac{\tilde f(T)}{(T+1)^{k+1}} \Big](z).
\end{align} 
Since $f$ and $\tilde f$ differ by a polynomial, we get \begin{align}\label{eqn:3.15}
\int_0^T \frac{ f(t)}{(t+1)^{k+1}} dt = \frac{C_{k+1}}{k!} \log(T+1)
+ L_1 + o(1)
\end{align} for some finite $L_1$. 

Next, we calculate the first integral of~\eqref{eqn:3.10}. 
We use the same strategy and apply Theorem~\ref{theorem:tauber1} to $\tilde \alpha_\pm(T) / (T+1)^k$. This gives \begin{align}\label{eqn:3.16}
\lim_{T\to \infty} \int_0^T \frac{\tilde \alpha_\pm(t)}{(t+1)^k} dt
= \lim_{z\to0} \Lcal \Big[ \frac{\tilde \alpha_\pm(T)}{(T+1)^k} \Big](z).
\end{align} 
Since $\alpha_\pm$ and $\tilde \alpha_\pm$ differ by polynomials, if $k\ge 1$, we have \begin{align}
\int_0^T \frac{ \alpha_+(t)}{(t+1)^k} dt 
&= \int_0^T  \frac{ A_{k+1} }{k!} \frac{ t^k}{ (t+1)^k} dt + \frac{ A_k}{(k-1)!} \log(T+1) +L_2 + o(1), \\
\int_0^T \frac{ \alpha_-(t)}{(t+1)^k} dt 
&= \int_0^T  \frac{ B_{k+1} }{k!} \frac{ t^k}{ (t+1)^k} dt + \frac{ B_k}{(k-1)!} \log(T+1) +L_3 + o(1), 
\end{align} for some finite $L_2$, $L_3$. 
We subtract the two equations and use $f' = \alpha_+ - \alpha_-$ and $A_{k+1} = B_{k+1}$ to get \begin{align}\label{eqn:3.19}
\int_0^T \frac{ f'(t)}{(t+1)^k} dt 
= \frac{ A_k - B_k}{(k-1)!} \log(T+1) +L_2 - L_3 + o(1).
\end{align}
Combining equations \eqref{eqn:3.10}, \eqref{eqn:3.15}, and \eqref{eqn:3.19}, we obtain \begin{align}
\frac{f(T)}{(T+1)^k}
= \frac{ A_k - B_k - C_{k+1} }{(k-1)!} \log(T+1) + L_4 + o(1)
\end{align}
for some finite $L_4$. 
Now since $f(T) = O(T^k)$, the left-hand side of the equation is bounded. This implies $A_k - B_k - C_{k+1} = 0$ and $\lim_{T\to\infty} \frac{f(T)}{(T+1)^k} = L_4$. 

If $k=0$, since $A_{k+1} = B_{k+1}$, we have $f' = \alpha_+ - \alpha_- = \tilde \alpha_+ - \tilde \alpha_-$. Denoting $\int_0^{\infty-} g = \lim_{T\to \infty} \int_0^T g$, equation~\eqref{eqn:3.16} gives \begin{align}
\int_0^{\infty-} f' 
= \int_0^{\infty-} \tilde \alpha_+ - \tilde \alpha_-
= \lim_{z\to 0} \Lcal [\tilde\alpha_+ - \tilde\alpha_- ](z)
= \lim_{z\to 0} \Lcal f'(z). 
\end{align}
As $z\to 0^+$, we have $\Lcal f'(z) = z\Lcal [f](z) - f(0) \to C-f(0)$. It then follows from $\int_0^{\infty-} f' = \lim_{T\to\infty} f(T) - f(0)$ that $\lim_{T\to\infty}f(T) = C$.

It remains to identify $L_4$ for the $k\ge 1$ case. Since $f\ge0$, the Hardy--Littlewood Tauberian theorem \cite[Theorem I.15.1]{Korevaar2004} states that $\Lcal f \sim Cz^{-(k+1)}$ as $z\to 0$ implies $\int_0^T f \sim \frac{C}{(k+1)!} T^{k+1}$ as $T\to \infty$. From this Ces\`aro sum and existence of $\lim_{T\to \infty} \frac{f(T) }{T^k}$, it is elementary to identify that $\lim_{T\to \infty} \frac{f(T) }{T^k} = \frac{C}{k!}$. 
\end{proof}

\subsection{Proof of Theorem~\ref{theorem:moments}}
We prove asymptotics for the numerator and the denominator of (\ref{eqn:6}) separately. The quotient of the limits then yields the theorem. We begin by proving two lemmas that will be used to verify the hypotheses of Theorem~\ref{theorem:tauber0}. 
The first lemma calculates derivatives, and it is proved in the same way as the Kolmogorov forward equations. The only difference is that we get an extra term for staying at site $j$.
\begin{lemma} \label{lemma:T-derivative}
Recall $P_{0j}(g, T) = E_0 [ e^{-g\sum_x \phi(L_{T,x})} \1_{X(T)=j} ]$. 
We have: 
\begin{enumerate}
\item[(i)] For any $j \in \Z$,
\begin{align} \label{eqn:P0j'}
\frac{\del P_{0j}}{\del T}
= -g E_0 \big[ \phi'(L_{T,j})e^{-g\sum_x \phi(L_{T,x})}\1_{X(T)=j} \big] 
+ ( P_{0, j-1} -2P_{0j} + P_{0, j+1}) .
\end{align}

\item[(ii)]
\begin{align}
\frac{\del}{\del T} \sum_{j=1}^\infty P_{0j}(g, T) =
&-g \sum_{j=1}^\infty E_0 \big[ \phi'(L_{T,j})e^{-g\sum_x \phi(L_{T,x})}\1_{X(T)=j} \big]  +( P_{00} - P_{01} ).
\end{align}

\item[(iii)]
For any $k\in\N$, 
\begin{multline}
\frac{\del}{\del T} \sum_{j=1}^\infty j^k P_{0j}(g, T) =
- g \sum_{j=1}^\infty j^k E_0 \big[ \phi'(L_{T,j})e^{-g\sum_x \phi(L_{T,x})}\1_{X(T)=j} \big] \\
+ P_{00} + (2^k-2)P_{01} + 2\sum_{j=2}^\infty  \bigg( \sum_{\substack{l=2 \\ l\ \text{even} }}^k  \binom{k}{l} j^{k-l}\bigg) P_{0j}. 
\end{multline}
\end{enumerate}
\end{lemma}

\begin{proof}
(i) Let $T\ge 0$ and $h>0$. Consider \begin{align}
P_{0j}(g, T+h) = E_0 [ e^{-g\sum_x \phi(L_{T+h,x})} \1_{X(T+h)=j} ]. \end{align}
We separate the right-hand side into three events. 
If the walk makes no jumps between time $T$ and $T+h$, then $X(T)=j$, and \begin{align}
L_{T+h,x} =\begin{cases}
L_{T,j}+h, & x=j, \\
L_{T,x}, &x\ne j. \end{cases}
\end{align} 
Since the jump rates of the walk are 1 to the left and 1 to the right, the probability for this event is $e^{-2h} = 1 -2h + o(h)$. 

If the walk makes exactly one jump between time $T$ and $T+h$, then $X(T) = j \pm 1$ with equal probability. For each of the starting points, the probability of jumping to $j$ by time $T+h$ is $1-e^{-h} = h + o(h)$. 
If the walk makes more than two jumps between time $T$ and $T+h$, this happens with $O(h^2)$ probability. Combining the three events, we get 
\begin{multline}
P_{0j}(g, T+h) = (1-2h)  E_0 [ e^{-g [ \phi(L_{T,j}+h)+\sum_{x\ne j} \phi(L_{T,x}) ]} \1_{X(T)=j} ]  \\
\quad +h P_{0,j-1}(g,T) + h P_{0,j+1}(g,T) + o(h).
\end{multline}
Thus, \begin{multline}
\lim_{h\to 0^+} \frac{P_{0j}(g, T+h) - P_{0j}(g, T)}{h}
= E_0 [ (-g)\phi'(L_{T,j})e^{-g \sum_x \phi(L_{T,x}) ]} \1_{X(T)=j} ]  \\
\quad - 2 P_{0j}(g,T) + P_{0,j-1}(g,T) + P_{0,j+1}(g,T).
\end{multline}
The left limit $h\to 0^-$ is proved similarly. 

\smallskip \noindent
(ii) This is from summing \eqref{eqn:P0j'} over $j \in \N$. 

\noindent
(iii) This is from summing $ j^k \cdot \eqref{eqn:P0j'}$. For $j\ge2$, using the Binomial Theorem, the coefficient for $P_{0j}$ is \begin{align}
(j-1)^k - 2j^k + (j+1)^k
= \sum_{l=1}^k \binom k l j^{k-l} [ (-1)^l + 1^l ]
= 2  \sum_{\substack{l=2 \\ l\text{ even} }}^k  \binom{k}{l} j^{k-l}, 
\end{align} which gives the desired result. 
\end{proof}

The next lemma establishes analyticity properties. The proof is algebraic. 
\begin{lemma} \label{lemma:chi_analytic} 
Let $g>0$ and $\nu_c = \nu_c(g)$. 
\begin{enumerate}
\item[(i)]
For any $0\ne y\in \R$, $1$ is not in the spectrum of $Q(g, \nu_c + iy)$. 

\item[(ii)]
The map $\nu \mapsto \chi_+(g, \nu) = \sum_{j=1}^\infty G_{0j}(g, \nu)$ can be extended to a meromorphic function on an open set containing the closed half-plane $\{ \nu \in \C \mid \Re(\nu)\ge \nu_c\}$, and it has a unique pole at $\nu=\nu_c$. 

\item[(iii)]
For any $k\in \N$, the map $\nu \mapsto \sum_{j=1}^\infty j^k G_{0j}(g, \nu)$ can be extended to a meromorphic function on an open set containing the closed half-plane  $\{ \nu \in \C \mid \Re(\nu)\ge \nu_c\}$, and it has a unique pole at $\nu=\nu_c$. 
\end{enumerate}
\end{lemma}

\begin{proof}
(i) Since $Q$ is compact, we only need to prove $1$ is not an eigenvalue. 
Denote $Q_c = Q(g, \nu_c)$ and $k_c(t, s) = k_0(t, s; g, \nu_c)$ for the kernel of $Q_c$. For $y\ne 0$, define an operator $Uf(t) = f(t)e^{-\half iyt}$. Then $Q = Q(g, \nu_c + iy)$ can be decomposed as \begin{align}
Qf(t) 
&= \int_0^\infty f(s)  e^{-\half iy t} e^{-\half iys} k_c(t,s)ds
= UQ_cUf(t). 
\end{align}

Suppose $Qf = (UQ_cU)f = 1f$, then by definition of $U$, we have 
$ e^{-\half iyt} Q_cUf(t) = f(t)$, so \begin{align}
 \langle Q_cUf, Uf \rangle = 
\int_0^\infty e^{\half iyt} f(t)  \overline{e^{-\half iyt} f(t)} dt
= \int_0^\infty e^{iyt} \abs{ f(t) }^2. \end{align}
Since $Q_c$ is a positive operator, we have \begin{align}
\langle Q_cU f, Q_cU f \rangle
\le \| Q_c\| \langle Q_cU f, U f \rangle 
= 1\cdot \int_0^\infty e^{iyt}\abs{ f(t) }^2 dt. \end{align}
But $\langle Q_cU f, Q_cU f \rangle =  \int_0^\infty e^{\half iyt} f(t)  \overline{e^{ \half iyt} f(t)} dt =  \int_0^\infty \abs{ f(t) }^2 dt, $
so plugging in and taking the real part give 
\begin{align}
0 \le \int_0^\infty (\cos(yt) - 1) \abs{ f(t) }^2 dt. \end{align}
For $y\ne0$, $\cos(yt) - 1<0$ for almost every $t$. This forces $f=0$ almost surely, so $1$ cannot be an eigenvalue of $Q$. 

\smallskip \noindent
(ii) By Proposition~\ref{prop:chi}(i), for $\Re(\nu) > \nu_c$, \begin{align}
\chi_+(g, \nu) 
= \left \langle Q(1-Q)\inv(q) , \bar q \right \rangle. \label{eqn:3.33}
\end{align}
By Proposition~\ref{prop:Gij}(i), $q$ is holomorphic in $\nu$ for $\Re(\nu) > \nu_c - \eps$. Since the conjugation $\bar q$ and the conjugation of the second argument of the inner product cancel each other, it suffices to prove that $(1-Q)\inv$ is well-defined on an open set containing the closed half-plane, except at $\nu=\nu_c$. 

By part (i), $(1-Q)\inv$ is well-defined at $\nu = \nu_c+iy$ with $y\ne 0$. By continuity, $(1-Q)\inv$ is well-defined on a small neighborhood around all such $\nu = \nu_c+iy$, $y\ne 0$. 
At $\nu=\nu_c$, we know 1 is the largest eigenvalue of $Q(g, \nu_c)$ and it is simple by Lemma~\ref{lemma:Q}(ii). Hence, the Kato--Rellich Theorem identifies the top of the spectrum $\lambda(\nu;g)$ of $Q(g, \nu)$ near $\nu=\nu_c$ (also see Remark~\ref{remark:kato}). Since for $\nu\in \R$ we have $\lambda(\nu;g) = \norm{Q(g, \nu)}$, the function $\lambda(\nu;g)$ must be non-constant. Thus, there exists a punctured neighborhood of $\nu=\nu_c$ in which $\lambda(\nu;g)\ne 1$. In this punctured neighborhood, $(1-Q)\inv$ is well-defined. 
Together, this proves that $(1-Q)\inv$ is well-defined on an open set containing the closed half-plane, except at $\nu=\nu_c$. 
The right-hand side of \eqref{eqn:3.33} then gives the desired extension. 

\smallskip
\noindent
(iii) For $k=1$, this follows from the same reasoning and equation~\eqref{eqn:2.54}, which states that \[
\sum_{j=1}^\infty j G_{0j}(g, \nu) 
=  \left \langle Q(1-Q)^{-2}(q) , \bar q \right \rangle, \qquad (\Re(\nu)>\nu_c). 
\] 
Similarly, for $k>1$, we use 
\begin{align}
\sum_{j=1}^\infty j^k G_{0j}(g, \nu) 
=  \left \langle \Big( \sum_{j=1}^\infty j^k Q^j \Big)(q) , \bar q \right \rangle, \qquad (\Re(\nu)>\nu_c). 
\end{align}
Note that $ \sum_{j=1}^\infty j^k Q^j$ belongs to the $\Z$-algebra generated by $Q$ and $(1-Q)\inv$ by holomorphic functional calculus, because $\sum_{j=1}^\infty j^k z^j$ belongs to the $\Z$-algebra generated by $z$ and $(1-z)\inv$ (by differentiating $ (1-z)\inv = \sum_{j=0}^\infty z^j$).
\end{proof}

\begin{proof}[Proof of Theorem~\ref{theorem:moments}]
We will see that $\theta(g) = 
\big(-\frac{\del \|Q\|}{\del\nu}\big\lvert_{\nu=\nu_c(g)}\big)\inv$. 
First, we handle the denominator of~\eqref{eqn:6}. 
Fix $g>0$. Define \begin{align}
f(T) =   \sum_{j=1}^\infty P_{0j}(g, T) e^{-\nu_c T} \ge 0, 
\end{align}
We differentiate $f$ using the product rule and Lemma~\ref{lemma:T-derivative}(ii). 
Let \begin{align}
 \alpha_+(T) &= -\nu_c \sum_{j=1}^\infty P_{0j}(g, T) e^{-\nu_c T}
+ P_{00}(g, T)e^{-\nu_c T}, \\
 \alpha_-(T) &= g \sum_{j=1}^\infty E_0 \big[ \phi'(L_{T,j})e^{-g\sum_x \phi(L_{T,x})}\1_{X(T)=j} \big] e^{-\nu_c T} + P_{01}(g, T)e^{-\nu_c T}, 
\end{align}
then $f' = \alpha_+ - \alpha_-$ and $ \alpha_\pm \ge 0$, because $\nu_c\le 0$ by Proposition~\ref{prop:chi}(ii) and $\phi'\ge 0$ by assumption~\eqref{A2}. 

Notice $\Lcal f (z) = \chi_+(g, \nu_c + z)$, so $\Lcal f(z)$ converges for $\Re(z)>0$ by Proposition~\ref{prop:chi} and extends to a meromorphic function on an open set containing the closed half-plane by Lemma~\ref{lemma:chi_analytic}. 
For $\Lcal \alpha_+(z) = -\nu_c \chi_+(g, \nu_c + z) + G_{00}(g, \nu_c+z)$, the same is true because $G_{00}(g, \nu_c+z)$ is holomorphic in $z$ for $\Re(z) > -\eps$ by Proposition~\ref{prop:Gij}. 
For $\Lcal \alpha_-$, using the method of Section~\ref{section:2pt} and assumption~\eqref{A4}, it is easy to prove that \begin{align}
 \Lcal E_0 \big[ \phi'(L_{T,j})e^{-g\sum_x \phi(L_{T,x})}\1_{X(T)=j} \big] = 
\left \langle Q^j [q](t), \phi'(t)  \overline{q(t)} \right\rangle. \end{align}
The required properties then follow analogously. 
We also know $\Lcal f$ has a unique simple pole at $z=0$, and each of 
$\Lcal \alpha_\pm$ either has a unique simple pole at $z=0$ or is holomorphic. 
Therefore, by Theorem~\ref{theorem:tauber0} with $k=0$ and the $\nu\to \nu_c^+$ asymptotics in \eqref{eqn:chi_asym}, we conclude \begin{align}
\sum_{j=1}^\infty P_{0j}(g, T) e^{-\nu_c T} \to 
\bar u \(-\frac{\del \|Q\|}{\del\nu}\Big\lvert_{\nu=\nu_c(g)}\)\inv, \qquad T\to \infty. \label{eqn:95}
\end{align}

For the numerator of~\eqref{eqn:6}, we use \begin{align}
f(T) 
&= \sum_{j=1}^\infty j^k P_{0j}(g, T) e^{-\nu_c T}, 
\\
\alpha_+(T)  
&= -\nu_c \sum_{j=1}^\infty j^kP_{0j}(g, T) e^{-\nu_c T}  \\
&\quad + \Big[ P_{00} + (2^k-2)P_{01} + 2\sum_{j=2}^\infty  \bigg( \sum_{\substack{l=2 \\ l\text{ even} }}^k  \binom{k}{l} j^{k-l}\bigg) P_{0j} \Big]
e^{-\nu_c T}, 
\nonumber \\
\alpha_-(T) 
&= g \sum_{j=1}^\infty j^k E_0 \big[ \phi'(L_{T,j})e^{-g\sum_x \phi(L_{T,x})}\1_{X(T)=j} \big] e^{-\nu_c T}.
\end{align}
This $\Lcal f$ has a pole of order $k+1$ by Corollary~\ref{corollary:correlation}. It follows from Theorem~\ref{theorem:tauber0} and the $\nu\to \nu_c^+$ asymptotics in \eqref{eqn:j^k_asymp} that \begin{align}
\sum_{j=1}^\infty j^k P_{0j}(g, T) e^{-\nu_c T} \sim 
\bar u \(-\frac{\del \|Q\|}{\del\nu}\Big\lvert_{\nu=\nu_c(g)}\)^{-(k+1)} T^k
, \qquad T\to \infty. 
\end{align}
Dividing by the denominator (\ref{eqn:95}), we get \[
\frac{ \sum_{j=1}^\infty j^k P_{0j}(g, T) }{ \sum_{j=1}^\infty P_{0j}(g, T) } \sim \(-\frac{\del \|Q\|}{\del\nu}\Big\lvert_{\nu=\nu_c(g)}\)^{-k} T^k, \qquad T\to \infty, \]
which is the desired result with $\theta(g) = 
\big(-\frac{\del \|Q\|}{\del\nu}\big\lvert_{\nu=\nu_c(g)}\big)\inv$. 
\end{proof}

\section{Monotonicity of speed} \label{section:monotonicity}
In this section, we prove that the escape speed $\theta(g) = \big(-\frac{\del \|Q\|}{\del\nu}\big\lvert_{\nu=\nu_c(g)}\big)^{-1}$ has a strictly positive derivative $\theta'(g)>0$. In Section~\ref{section:monotonicity-reduction}, we first build a sequence $\{ c_n \}_{n=0}^\infty$ for every $g$, and we prove that $\theta'(g)>0$ is equivalent to $c_0 + 2 \sum_{n=1}^\infty c_n >0$. Here, $c_n$ is related to the $n$-th power of the operator $Q$ evaluated at $(g, \nu_c(g))$ (see \eqref{def:cn}). 
In Section~\ref{section:monotonicity-dominance}, we use stochastic dominance to prove $c_0>0$ and $c_n\ge 0$ for all $n$, which imply $\theta'(g)>0$ and complete the proof of Theorem~\ref{theorem:WLLN}. 

Since we only encounter the operator $Q$ in this section, for simplicity, we write $k(t,s) = k_0(t,s)$ for the kernel of $Q$ (see \eqref{def:Q}). Also, we assume $\nu\in \R$ throughout the section.

\subsection{Reduction} \label{section:monotonicity-reduction}
We begin by calculating $\theta'(g)$ using the Implicit Function Theorem. Using subscripts for partial derivatives and denoting $\lambda = \norm{Q}$, we have 
\begin{align} \label{eqn:101}
\frac{d}{dg}\(  \frac{1}{\theta(g)}  \) 
&= -\lambda_{\nu g} - \lambda_{\nu\nu}\frac{d\nu_c}{dg}  \\
&= -\lambda_{\nu g} - \lambda_{\nu\nu}(-\lambda_g / \lambda_\nu)  \nonumber \\
&= (  \lambda_{\nu g}\lambda_\nu  - \lambda_{\nu\nu}\lambda_g)/(-\lambda_\nu).\nonumber 
\end{align}
Since $\lambda_\nu <0$ by Lemma~\ref{lemma:nu_c}, we have $\theta' >0$ if and only if $\lambda_{\nu g}\lambda_\nu  - \lambda_{\nu\nu}\lambda_g <0$. 
This combination of derivatives is central to the reduction so we give it a name. For a $C^2$ function $F = F(g, \nu)$, we define 
\begin{align}\label{def:L}
L[F] =  F_{\nu g}F_\nu  - F_{\nu\nu}F_g. \end{align}
The goal is to prove $L[\lambda] < 0$. 
However, it is difficult to calculate $L[\lambda]$ directly because of the second derivatives. Instead, as suggested in \cite{GH1993}, we can calculate $L[ \langle Q^n(g, \nu) f, f \rangle^{1/n} ]$ for some $f>0$ and then send $n\to\infty$. The idea is to utilize $\|Q(g, \nu) \| = 
\lim_{n\to\infty}\langle Q^n(g, \nu) f, f \rangle^{1/n}$. This is justified by the following lemma, with $f$ chosen to be a leading eigenvector of $Q$. 

\begin{lemma} \label{lemma:Hn}
Fix $g_0>0$ and $\nu_0 = \nu_c(g_0)$. Let $h_0$ be the positive normalized leading eigenvector of $Q(g_0, \nu_0)$, \ie, $h_0$ satisfies $Q(g_0, \nu_0)h_0 = h_0$, $h_0>0$ and $\| h_0\|_2 = 1$. For any $n\ge1$, define \begin{align} \label{def:Hn}
H_n(g, \nu) = \langle Q^n(g, \nu)h_0, h_0 \rangle^{1/n}. \end{align}
Then for $*=g, \nu$, we have $ \del_* H_{n} (g_0, \nu_0) = \lambda_*$ and $\del_*\del_\nu H_{n}(g_0, \nu_0) \to \lambda_{\nu *}$ as $n\to \infty$. 
Hence, \begin{align}\label{LH_limit}
\lim_{n\to \infty} L[H_n] \big \lvert_{g_0, \nu_0} = L[\lambda] \big\lvert_{g_0, \nu_0}. \end{align}
\end{lemma}
This lemma is proved by calculating $L[H_n]$ using differentiation rules and calculating $L[\lambda]$ using Cauchy's integral formula. The proof is given in Appendix~\ref{appendix:Hn}. 
The function $H_n(g, \nu)$ is more tractable than the operator norm $\lambda(g, \nu) = \norm{Q(g, \nu)}$. 

\begin{lemma} \label{lemma:L[Hn]}
Assume the hypotheses of Lemma~\ref{lemma:Hn}. 
For any $n\ge 0$, define 
\begin{align} \label{def:cn}
c_n &= 
\Big\langle Q^{n}\big[t  h_0(t)\big] (s), \phi(s)  h_0(s)\Big\rangle
\cdot \int_0^\infty s h_0^2(s) ds  \\
&\quad- \Big\langle Q^{n}\big[t  h_0(t)\big](s), s  h_0(s)\Big\rangle
\cdot \int_0^\infty \phi(s) h_0^2(s) ds, \nonumber
\end{align}
with $Q$ evaluated at $(g_0, \nu_0)$. 
Then for all $n\ge 1$, \begin{align}
L[H_n] \big \lvert_{g_0, \nu_0} 
&= -\frac 1 n 
\bigg( - \half c_0 + \half c_n + \sum_{i=1}^n \sum_{j=1}^n c_{\abs{ j-i }} \bigg).
\end{align}
\end{lemma}

Combining the two lemmas, we sum diagonally (the convergence of the series is controlled by the second largest eigenvalue, which is $<1$) to obtain \begin{align}
L[\lambda] \big\lvert_{g_0, \nu_0}
= \lim_{n\to\infty} L[H_n] \big \lvert_{g_0, \nu_0} 
= -c_0 - 2\sum_{n=1}^\infty c_n. 
\end{align}
In the next subsection, we will prove $c_0>0$ and $c_n\ge 0$ for all $n$. This will allow us to conclude $L[\lambda] \big\lvert_{g_0, \nu_0}<0$, which is equivalent to $\theta'(g_0)>0$. 

The proof of Lemma~\ref{lemma:L[Hn]} uses the following lemma due to \cite{GH1993}.

\begin{lemma} \label{lemma:phi}
Let $F=F(g, \nu)$ be $C^2$ and $\phi:\R\to\R$ be differentiable on the image of $F$. Then \begin{align}
L[\phi(F)] = (\phi'(F))^2 L[F]. \end{align}
\end{lemma}

\begin{proof}
This is a direct calculation using differentiation rules. We have
\begin{align*}
L[\phi(F)]
&= (\phi(F))_{\nu g}(\phi(F)) _\nu - (\phi(F))_{\nu\nu}(\phi(F))_g   \\
&= (\phi'(F)F_\nu)_g(\phi'(F)F_\nu) - (\phi'(F)F_\nu)_\nu (\phi'(F)F_g)\nonumber \\
&= ({\phi'' (F) F_gF_\nu} +\phi' (F) F_{\nu g})(\phi' (F) F_\nu) 
- ({\phi'' (F) F_\nu F_\nu} +\phi' (F) F_{\nu \nu})(\phi' (F) F_g) \nonumber \\
&= (\phi'(F))^2 (F_{\nu g}F_\nu - F_{\nu\nu}F_g) \nonumber \\
&= (\phi'(F))^2 L[F], \nonumber
\end{align*} which is the desired result. 
\end{proof}

\begin{proof}[Proof of Lemma~\ref{lemma:L[Hn]}]
We use Lemma~\ref{lemma:phi} with $\phi(t) = t^{1/n}$ and $F(g, \nu) = \langle Q^n(g, \nu) h_0, h_0 \rangle$. By the hypotheses, we have $\langle Q^n(g_0, \nu_0) h_0, h_0 \rangle = \langle h_0, h_0 \rangle = 1$, hence, \begin{align} \label{eqn:4.11}
L[H_n] \big \lvert_{g_0, \nu_0} 
&= 
\( \frac{1}{n}(1) \)^2 L\big[ \langle Q^n(g, \nu) h_0, h_0 \rangle \big] \Big\lvert_{g_0, \nu_0} .
\end{align}
We calculate the right-hand side next. 
Since $h_0$ is fixed, we only need to differentiate $Q$. 
Recall $Q$ was defined in \eqref{def:Q} by \begin{align*}
Qf(t) &= \int_0^\infty f(s) k(t,s) ds, \\
k(t,s) &= \sqrt{p(t)}\sqrt{p(s)} e^{-t}e^{-s} I_0(2\sqrt{st}), \\
\sqrt{p(t)} &= e^{-\half g\phi(t) - \half \nu t}. 
\end{align*}
Observe $g$ and $\nu$ enter the kernel only through $\sqrt{p}$. 
Writing out all $n$ integrals in $Q^n h_0$ yields \begin{align} \label{eqn:Qn}
\langle Q^n  h_0, h_0 \rangle
&= \int_{(0, \infty)^{n+1}} h_0(s_n) k(s_n, s_{n-1}) \cdots k(s_1, s_0) \cdot  h_0(s_0) d\bf s,
\end{align} 
where $d\bf s = ds_0 \cdots ds_n$. 
Thus, the $g$-derivative is  
\begin{multline}
\frac{\del}{\del g}\langle Q^n  h_0, h_0 \rangle
= - \int_{(0, \infty)^{n+1}}  \bigg( \half \phi(s_n) + \sum_{j=1}^{n-1} \phi(s_j) + \half \phi(s_0)\bigg)
 \\
\cdot  h_0(s_n) k(s_n, s_{n-1}) \cdots k(s_1, s_0) \cdot  h_0(s_0) d\bf s.
\end{multline} 
When evaluated at $(g_0, \nu_0)$, we have $Q(g_0,\nu_0)h_0 = h_0$, so for each $j$,  \begin{align}
\int_{(0, \infty)^{n+1}}   \phi(s_j) h_0(s_n) k(s_n, s_{n-1}) \cdots k(s_1, s_0) \cdot  h_0(s_0) d\bf s
= \int_0^\infty \phi(s_j) h_0^2(s_j) ds_j. \end{align}
This gives 
\begin{align}
\frac{\del}{\del g}\langle Q^n  h_0, h_0 \rangle \bigg\lvert_{g_0, \nu_0}
&= -n \int_0^\infty \phi(s) h_0^2(s) ds. 
\end{align} 

The $\nu$-derivative is similar and gives a multiplier of 
$ -\big( \half s_n + \sum_{j=1}^{n-1} s_j + \half s_0\big)$ to the integrand of  \eqref{eqn:Qn}. We get \begin{align}
\frac{\del}{\del \nu}\langle Q^n  h_0, h_0 \rangle \bigg\lvert_{g_0, \nu_0}
&= -n \int_0^\infty s h_0^2(s) ds. \end{align}

For the second derivatives, we define \begin{align} \label{def:alpha}
\alpha_j =\alpha_j(n)= \begin{cases} 
\half &j=0, \\ 
1& 0<j<n, \\
\half & j=n,
\end{cases}\end{align}
then \begin{multline}
\frac{\del^2}{\del \nu \del g}\langle Q^n  h_0, h_0 \rangle \\
= 
\sum_{i=0}^n \sum_{j=0}^n \alpha_i \alpha_j 
\int_{(0, \infty)^{n+1}} \phi(s_i) s_j h_0(s_n) k(s_n, s_{n-1}) \cdots k(s_1, s_0) \cdot  h_0(s_0) d\bf s.
\end{multline}
When evaluating at $(g_0, \nu_0)$, we have $Q(g_0,\nu_0)h_0 = h_0$, so 
\begin{align}
&\int_{(0, \infty)^{n+1}} \phi(s_i) s_j h_0(s_n) k(s_n, s_{n-1}) \cdots k(s_1, s_0) \cdot  h_0(s_0) d\bf s  \\
&\qquad= \Big\langle Q^{\abs{ j-i }}\big[s_j h_0(s_j)\big](s_i), \phi(s_i) h_0(s_i)\Big\rangle. \nonumber \end{align}
Thus, \begin{align}
\frac{\del^2}{\del \nu \del g}\big\langle Q^n  h_0, h_0 \big\rangle \bigg\lvert_{g_0, \nu_0}&= 
\sum_{i=0}^n \sum_{j=0}^n \alpha_i \alpha_j 
\Big\langle Q^{\abs{ j-i }}\big[s_j h_0(s_j)\big](s_i), \phi(s_i) h_0(s_i)\Big\rangle. \end{align}
Similarly, \begin{align}
\frac{\del^2}{\del \nu^2 }\big\langle Q^n  h_0, h_0 \big\rangle \bigg\lvert_{g_0, \nu_0}&= 
\sum_{i=0}^n \sum_{j=0}^n \alpha_i \alpha_j 
\Big\langle Q^{\abs{ j-i }}\big[s_j h_0(s_j)\big](s_i), s_i h_0(s_i)\Big\rangle. \end{align}

Using the definition of $c_n$ in~\eqref{def:cn} and the definition of $L$ in~\eqref{def:L},  combining the sums above gives 
\begin{align} 
L\big[ \langle Q^n(g, \nu) h_0, h_0 \rangle \big] \Big\lvert_{g_0, \nu_0} 
&= \sum_{i=0}^n \sum_{j=0}^n \alpha_i \alpha_j (-n) c_{\abs{ j-i }} \\
&= (-n) \bigg( -\half c_0 +\half c_n + 
\sum_{i=1}^n \sum_{j=1}^n  c_{\abs{ j-i }} \bigg). \nonumber
\end{align} 
This and equation~\eqref{eqn:4.11} give the claim. 
\end{proof}

\subsection{Stochastic dominance} \label{section:monotonicity-dominance}
In this section, we prove $c_0>0$ and $c_n\ge 0$, which then imply $\theta'(g_0)>0$ by Section~\ref{section:monotonicity-reduction}. 
In the following proposition, we rewrite the inequality $c_n\ge 0$ in a quotient form (recall the definition of $c_n$ in \eqref{def:cn}). This allows the inequality to be interpreted as an inequality between the expectations of two random variables.

\begin{proposition} \label{prop:Qn-inequality}
Fix $g_0>0$ and $\nu_0 = \nu_c(g_0)$. Let $Q = Q(g_0, \nu_0)$. Let $h_0$ be the positive normalized leading eigenvector of $Q$, \ie, $h_0$ satisfies $Qh_0 = h_0$, $h_0>0$ and $\| h_0\|_2 = 1$. Then for any $n\ge 0$, we have \begin{equation}\label{eqn:Qn-inequality}
\frac{\Big\langle Q^n\big[t h_0(t)\big](s), \phi(s) h_0(s)\Big\rangle}{\displaystyle \int_0^\infty \phi(s) h_0^2(s) ds} 
\ge
\frac{\Big\langle Q^n\big[t  h_0(t)\big](s), s h_0(s)\Big\rangle}{\displaystyle \int_0^\infty s h_0^2(s) ds},
\end{equation}
and the inequality is strict for $n=0$. 
\end{proposition}
To illustrate the method of the proof, we will first prove the case $n=0$ using (first-order) stochastic dominance. For real-valued random variables $X, Y$, we write $X\stle Y$ if $P(X > x) \le P(Y>x)$ for all $x\in \R$. If $X, Y$ have density functions $f_X, f_Y$ respectively, a sufficient condition for $X\stle Y$ is that $f_Y/f_X$ is an increasing function. A consequence of $X\stle Y$ is that $EX \le EY$. 

\begin{proof}[Proof for the case $n=0$]
If $n=0$, the goal (\ref{eqn:Qn-inequality}) reduces to \[
\frac{ \int_0^\infty s\phi(s) h_0^2(s) ds}{ \int_0^\infty \phi(s) h_0^2(s) ds} 
>
\frac{ \int_0^\infty s\cdot s h_0^2(s) ds} { \int_0^\infty s h_0^2(s) ds}. \]
This can be written as $EY > EX$, 
where $Y, X$ are random variables on $(0, \infty)$ defined by the density functions 
\begin{align}
f_Y(s) = \frac{\phi(s) h_0^2(s)}{\int_0^\infty \phi(s) h_0^2(s) ds}, \qquad
f_X(s) = \frac{s h_0^2(s)}{\int_0^\infty s h_0^2(s) ds}. 
\end{align}
Notice $f_Y(s)/f_X(s) = c\phi(s)/s$ for some positive constant $c$, so it is increasing by assumption~\eqref{A2}. Thus, $X\stle Y$ and $EX \le EY$. The strict inequality follows from $X$ and $Y$ having different distributions, which is a consequence of $h_0(s)>0$ for all $s\ge0$. 
\end{proof}

To prove the general case, we need a result on multivariate stochastic order. For random vectors $\bf X, \bf Y \in \R^{n+1}$, we say $\bf X \stle \bf Y$ if $P(\bf X\in U) \le P( \bf Y \in U )$ for any increasing set $U\subset \R^{n+1}. $
\begin{lemma}[{\cite[Theorem 6.B.3]{Shaked-Shanthikumar2007}}] \label{thm:multi}
Let $\bf X = (X_0, \dots, X_n)$ and $\bf Y = (Y_0, \dots, Y_n)$ be $\R^{n+1}$-valued random variables. If $X_0 \stle Y_0$ and the conditional distributions satisfy \begin{equation}\label{conditional_dominance}
[X_i\mid X_0 = x_0, \dots, X_{i-1} = x_{i-1} ] 
\stle [Y_i\mid Y_0 = y_0, \dots, Y_{i-1} = y_{i-1} ] \end{equation}
whenever $(x_0, \dots, x_{i-1}) \le (y_0, \dots, y_{i-1})$ for all $i= 1, \dots, n$, 
then $\bf X \stle \bf Y$. 
As a result, $X_n\stle Y_n$ and $EX_n \le EY_n$. 
\end{lemma}

\begin{proof}[Proof of Proposition~\ref{prop:Qn-inequality} for the case $n>0$]
Recall $Q$ was defined by
$Qf(t) = \int_0^\infty f(s) k_0(t,s) ds$ and we write $k(t,s) = k_0(t,s)$. 
With $s_n$ replacing $t$ and $s_0$ replacing $s$, the numerator of the left-hand side of the goal \eqref{eqn:Qn-inequality} is 
\begin{multline}
\Big\langle Q^n\big[t \cdot h_0(t)\big](s), \phi(s) h_0(s)\Big\rangle \\
= \int_{(0, \infty)^{n+1}} s_n\cdot h_0(s_n) k(s_n, s_{n-1}) \cdots k(s_1, s_0) \cdot \phi(s_0) h_0(s_0) d\bf s. 
\end{multline}
We define a random variable $\bf Y = (Y_0, \dots, Y_n) \in \R_+^{n+1}$ by the density function \begin{align}
f_{\bf Y}(s_0, \dots, s_n) = \frac{1}{Z_Y} h_0(s_n) k(s_n, s_{n-1}) \cdots k(s_1, s_0) \cdot \phi(s_0) h_0(s_0). 
\end{align} 
Since  the normalizing constant $Z_Y$ is given by
\begin{align}
Z_Y &= \int_{(0, \infty)^{n+1}}  h_0(s_n) k(s_n, s_{n-1}) \cdots k(s_1, s_0) \cdot \phi(s_0) h_0(s_0)  \\
&= \langle Q^n h_0, \phi h_0 \rangle
= \langle  h_0, \phi h_0 \rangle
= \int_0^\infty \phi(s_0) h_0^2(s_0)ds_0, \nonumber
\end{align} the left-hand side of \eqref{eqn:Qn-inequality} is exactly $EY_n$. 
Similarly, we define a random variable $\bf X = (X_0, \dots, X_n ) \in \R_+^{n+1}$ by the density function \begin{align}
f_{\bf X}(s_0, \dots, s_n) = \frac{1}{Z_X} h_0(s_n) k(s_n, s_{n-1}) \cdots k(s_1, s_0) \cdot s_0 h_0(s_0), 
\end{align} then the right-hand side of \eqref{eqn:Qn-inequality} is exactly $EX_n$.
We use Lemma~\ref{thm:multi} to prove $EY_n \ge EX_n$. 

To see that $X_0 \stle Y_0$, notice the density functions for $X_0$ and $Y_0$ are \begin{align}
f_{X_0}(s_0) = \frac{1}{Z_X} s_0 h_0^2(s_0), \qquad
f_{Y_0}(s_0) = \frac{1}{Z_Y} \phi(s_0) h_0^2(s_0) 
\end{align}
because $Qh_0=h_0$. The dominance is proved in the $n=0$ case. 
To show the dominance of conditional distributions (\ref{conditional_dominance}), we use conditional density functions. The conditional density of $X_i$ given $ X_0 = x_0, \dots, X_{i-1} = x_{i-1} $ is \begin{align}
f_i(s_i\mid x_0, \dots, x_{i-1}) = \frac{1}{Z_1}h_0(s_i) k(s_i, x_{i-1}) \cdots k(x_1, x_0)\cdot x_0 h_0(x_0),
\end{align}
and the conditional density function of $Y_i$ given $Y_0 = y_0, \dots, Y_{i-1} = y_{i-1}$ is \begin{align}
g_i(s_i\mid y_0, \dots, y_{i-1}) = \frac{1}{Z_2}h_0(s_i) k(s_i, y_{i-1}) \cdots k(y_1, y_0)\cdot \phi(y_0) h_0(y_0). 
\end{align}
By definition of $k(\cdot, \cdot) = k_0(\cdot,\cdot)$ in \eqref{def:k_j}, we have 
\begin{align}
\frac{g_i(s_i\mid y_0, \dots, y_{i-1})}{f_i(s_i\mid x_0, \dots, x_{i-1})}
&= C_1 \frac{k(s_i, y_{i-1})}{k(s_i, x_{i-1})}  \\
&= C_1 \frac{\sqrt{p(y_{i-1})}\sqrt{p(s_i)} e^{-y_{i-1}}e^{-s_i} I_0(2\sqrt{s_i y_{i-1}})}
{\sqrt{p(x_{i-1})}\sqrt{p(s_i)} e^{-x_{i-1}}e^{-s_i} I_0(2\sqrt{s_i x_{i-1}})} \nonumber \\
&= C_2 \frac{I_0(2\sqrt{s_i y_{i-1}})}{I_0(2\sqrt{s_i x_{i-1}})}, \nonumber
\end{align} where $C_1, C_2>0$ are constants not depending on $s_i$. 
To prove the required stochastic dominance, we show that this ratio is increasing in $s_i$. To simplify the notation, we prove \begin{align}
F(t) = \frac{I_0(\lambda_2 t)}{I_0(\lambda_1 t)} \end{align} to be increasing in $t\ge0$ if $0\le\lambda_1\le \lambda_2$, then we can apply this with $t = \sqrt{s_i}$, $\lambda_1 = 2\sqrt{x_{i-1}}$, and $\lambda_2 = 2\sqrt{y_{i-1}}$. Notice \begin{align}
(\log F(t))' = \lambda_2 \frac{I_0'(\lambda_2 t)}{I_0(\lambda_2 t)} - \lambda_1 \frac{I_0'(\lambda_1 t)}{I_0(\lambda_1 t)}. \end{align}
Since $\lambda_2\ge \lambda_1\ge 0$, $I_0' = I_1$, and since ${I_1}/{I_0}$ is non-negative and increasing\footnote{The monotonicity of the ratio $I_1/I_0$ was proved in, \eg,  \cite[Theorem 1, 1(b)]{segura2021monotonicity}.}, we get $(\log F(t))' \ge 0$ for $t\ge0$. Since $F>0$, we get $F'\ge 0$ as well. This verifies all hypotheses of Lemma~\ref{thm:multi} and finishes the proof. 
\end{proof}

\section*{Acknowledgements}
This work was supported in part by NSERC of Canada. 
We thank Gordon Slade for discussions and very helpful advice. 
We thank Roland Bauerschmidt for introducing the problem and the transfer matrix approach to us. 
We thank Frank den Hollander for encouragement and comments.

\begin{appendices}
\section{Supersymmetric representation of random walk}
\label{appendix:supersymmetry}
What we call the supersymmetric representation of random walk is an integral representation for certain functionals of local times of a continuous-time random walk. It is also known as the supersymmetric version of the BFS--Dynkin isomorphism theorem \cite[Corollary~11.3.7]{BBS2019}. 
The integral involved is an integration of differential forms. 
Further background of the isomorphism theorem can be found in \cite[Chapter~11]{BBS2019}.

\subsection{Differential forms, bosons, and fermions}
\subsubsection{Integration of differential forms}
Let $\Lambda = \{ 1, \dots, \abs{ \Lambda }\}$ be a finite set, which will be the state space of the random walk. We consider 2-component real fields over $\Lambda$. For each $x\in \Lambda$, let $(u_x, v_x)$ be real coordinates. The 1-forms $\{ du_x, dv_x \}_{x\in \Lambda}$ generate the Grassmann algebra of differential forms on $\R^{2\Lambda}$, with multiplication given by the anti-commuting wedge product. 

We write $u = (u_x)_{x\in \Lambda}, v = (v_x)_{x\in \Lambda}$. For $p\in \N_0$, a $p$-\emph{form} is a function of $(u, v)$ multiplied with a product of $p$ differentials. A \emph{form} $K$ is a sum of $p$-forms with possibly different values of $p$. The largest such $p$ is called the \emph{degree} of $K$. For a fixed $p$, the contribution of $p$-forms to $K$ is called the \emph{degree-}$p$ part of $K$. A form is called \emph{even} if it is the sum of $p$-forms with even $p$ only. Due to anti-commutativity, we can always write the degree-$2\lvert \Lambda \rvert$ part of a form $K$ as \begin{align}
 f(u,v) du_1\wedge dv_1 \wedge \cdots \wedge du_{\abs{ \Lambda }} \wedge dv_{\abs{ \Lambda }}. 
\end{align}
The integral of $K$ is then defined to be the Lebesgue integral
\begin{align}
\int_{\R^{2\lvert \Lambda \rvert}} K = \int_{\R^{2\lvert \Lambda \rvert}} f(u,v) du_1 dv_1  \cdots du_{\abs{ \Lambda }}  dv_{\abs{ \Lambda }}.
\end{align} 
Notice if the degree of $K$ is strictly less than $2\lvert \Lambda \rvert$, then its degree-$2\lvert \Lambda \rvert$ part is zero, so its integral is zero.

\subsubsection{Bosons and Fermions}
It is convenient to use complex coordinates. For $x\in \Lambda$, we define 
\begin{alignat}2
\phi_x &= u_x + i v_x, \qquad\quad~~ \bar\phi_x &&= u_x - i v_x, \\
d\phi_x &= du_x + i dv_x, \qquad  d\bar\phi_x &&= du_x - i dv_x. \nonumber
\end{alignat}
We call $(\phi, \bar \phi) = (\phi_x, \bar \phi_x)_{x\in \Lambda}$ the \emph{boson field}. 
We also define 
\begin{align}
\psi_x = \frac{1}{\sqrt{2\pi}} e^{-i\pi/4} d\phi_x, \qquad
\bar \psi_x = \frac{1}{\sqrt{2\pi}} e^{-i\pi/4}  d\bar \phi_x,
\end{align}
and call $(\psi, \bar \psi) = (\psi_x, \bar \psi_x)_{x\in \Lambda}$ the \emph{fermion field}. Then \begin{align}
\bar \psi_x \wedge \psi_x
= \frac 1 {2\pi i} d\bar \phi_x \wedge d\phi_x
= \frac 1 \pi du_x \wedge dv_x. 
\end{align}
The combination $\Phi = (\phi_x, \bar \phi_x, \psi, \bar \psi_x)_{x\in \Lambda}$ is called a \emph{superfield}. 
From now on, we drop the wedge symbol in the wedge product. 
One important field of forms is \begin{align}
\Phi^2 = (\Phi_x^2)_{x\in \Lambda} = ( \phi_x \bar \phi_x + \psi_x\bar \psi_x )_{x\in \Lambda}. \end{align}
For a complex $\abs{ \Lambda }\times \abs{ \Lambda }$ matrix $\Delta$, we define \begin{align}
( \Phi, -\Delta \Phi) = \sum_{x\in \Lambda} \( \phi_x (-\Delta \bar \phi )_x +  \psi_x (-\Delta \bar \psi )_x \). 
\end{align}

\subsubsection{Function of forms}
For $p\in \N$, consider a $C^\infty$ function $F:\R^p \to \R$. Let $K = (K_1, \dots, K_p)$ be a collection of even forms. Assume the degree-0 part $K_j^0$ of each $K_j$ is real. Then we define the form $F(K)$ using the Taylor series about the degree-0 part $K^0$, as 
\begin{align}
F(K) = \sum_\alpha \frac {1} {\alpha!} F^{(\alpha)}(K^0) (K-K^0)^\alpha, 
\end{align}
where the sum is over all multi-indexes $\alpha = (\alpha_1, \dots, \alpha_p)$, and $\alpha! = \prod_{j=1}^p \alpha_j!$, $(K-K^0)^\alpha = \prod_{j=1}^p (K_j - K_j^0)^{\alpha_j}$ (the order of the product does not matter because all $K_j$ are even). The sum is always finite because $(K_j - K_j^0)^{\alpha_j} = 0$ for all $\alpha_j > 2\lvert \Lambda \rvert$ by anti-commutativity. 
The key example is the following. With $p=1$ and $x\in \Lambda$, 
\begin{align} \label{eqn:A9}
F(\Phi_x^2) = F( \phi_x \bar \phi_x + \psi_x\bar \psi_x)
= F( \phi_x \bar \phi_x) + F'( \phi_x \bar \phi_x) \psi_x\bar \psi_x.
\end{align}

\subsection{Isomorphism theorem and supersymmetry}
Let $\{ X(t) \}_{t\ge 0}$ be a continuous-time random walk on a finite set $\Lambda$ with generator $\Delta$. 
We denote its expectation by $E_i$ if $X(0)=i$. 
The local time of $X$ at $x\in \Lambda$ up to time $T$ is defined by \begin{equation}
L_{T,x} = \int_0^T \1_{X(s)=x}\ ds. \end{equation}
We write $L_T = (L_{T,x})_{x\in \Lambda}$. 
The supersymmetric BFS--Dynkin isomorphism theorem \cite[Corollary 11.3.7]{BBS2019} relates local times of $X(t)$ with boson and fermion fields, as follows. 
\begin{theorem}[BFS--Dynkin isomorphism theorem] \label{theorem:isomorphism}
Let $F:\R^{\abs{ \Lambda }} \to \R$ be such that $e^{\eps \sum_{x\in \Lambda} t_x} F(t)$ is a Schwartz function for some $\eps>0$. Then \begin{align}
\int_0^\infty E_i (F(L_T) \1_{X(T)=j })\ dT
= \int_{\R^{2\lvert \Lambda \rvert}}\bar \phi_i \phi_j e^{-(\Phi, -\Delta \Phi)} F(\Phi^2), \label{eqn:isomorphism}
\end{align}
where $\{ X(t) \}_{t\ge 0}$ is a continuous-time random walk on $\Lambda$ with generator $\Delta$. 
\end{theorem}
We will use this theorem on the finite-volume two-point function $G_{ij}^N$ defined in \eqref{def:Gij^N}. Choosing the nearest-neighbor Laplacian $\Delta$ on the right-hand side allows us to use the transfer matrix approach. 

There is a symmetry between bosons and fermions called \emph{supersymmetry}. The next theorem is a demonstration of this. Notice the form $\Phi^2_x = \phi_x \bar \phi_x + \psi_x \bar \psi_x$ is unchanged if we interchange $(\phi_x, \bar \phi_x)$ with $(\psi_x, \bar \psi_x)$, so the integrands of the two sides of \eqref{eqn:A.12} are related by an interchange of bosons and fermions. 
For general results and discussions on supersymmetry, we refer to \cite[Section 11.4]{BBS2019}.

\begin{theorem} \label{theorem:supersymmetry}
Let $x\in \Lambda$ and $F: [0, \infty) \to \R$ be smooth. If $\lim_{t\to\infty} tF(t) = 0$ and the integrals exist, then \begin{align} \label{eqn:A.12}
\int_{\R^2} \bar \phi_x \phi_x F(\Phi_x^2) 
=\int_{\R^2} \bar \psi_x \psi_x F(\Phi_x^2) . 
\end{align}
\end{theorem}

\begin{proof}
Since bosons commute and fermions anti-commute, \begin{align}
\bar \phi_x \phi_x - \bar \psi_x \psi_x
= \phi_x  \bar \phi_x + \psi_x \bar \psi_x
= \Phi_x^2, 
\end{align}
so it is sufficient to prove $\int_{\R^2 }  \Phi_x^2 F(\Phi_x^2) =0$. By definition of the integral and by \eqref{eqn:A9}, \begin{align}
\int_{\R^2 }  \Phi_x^2 F(\Phi_x^2) 
&= \int_{\R^2 } (  \phi_x  \bar \phi_x + \psi_x \bar \psi_x ) ( F( \phi_x \bar \phi_x) + F'( \phi_x \bar \phi_x) \psi_x\bar \psi_x ) \\
&= \int_{\R^2 } [ F( \phi_x \bar \phi_x) + \phi_x  \bar \phi_x F'( \phi_x \bar \phi_x) ] \psi_x\bar \psi_x \nonumber \\
&= \int_{\R^2 } [ F( u^2 + v^2) + (u^2 + v^2) F'( u^2 + v^2) ] \frac{-1}{\pi} dudv \nonumber \\
&= - \int_0^\infty [ F( r^2) + r^2 F'( r^2) ]  dr^2 \nonumber \\
&= - \int_0^\infty \frac{d}{dt}\big( tF(t) \big)  dt = 0, \nonumber
\end{align}
which is the desired result. 
\end{proof}

\subsection{Proof of Proposition~\ref{prop:Gij^N}} \label{appendix:Gij^N}
We first prove Proposition~\ref{prop:Gij^N} for $\nu\in\R$ using the supersymmetric representation. We need the following lemma, which is a corollary of Proposition 2.5 and Lemma 2.6 in \cite{BS2020}. 
Since we need to deal with two superfields at the same time, we use the notation $D\Phi$ to signify that the integration is with respect to the superfield $\Phi$. We also recall the operators $T$ and $Q$ defined in~\eqref{def:T} and \eqref{def:Q}. 
Notice $f(t) = \sqrt{p(t)} = e^{-\half g\phi(t) - \half \nu t}$ satisfies the hypotheses of the lemma, because $\phi(t) \ge 0$ and $\phi(0)=0$ by assumption~\eqref{A0}. 

\begin{lemma} \label{lemma:outside_integral}
Let $\nu\in \R$. Fix a superfield $Z = (\zeta, \bar \zeta, \xi, \bar \xi)$, where $\xi = \frac 1 {\sqrt{2\pi }} e^{-i\pi/4} d\zeta$ and  $\bar \xi = \frac 1 {\sqrt{2\pi}} e^{-i\pi/4}  d\bar \zeta$. 
Let $f:[0, \infty)\to [0, \infty)$ be a smooth function such that $\sqrt{p} \cdot f$ is bounded. Then:

\begin{enumerate}
\item[(i)]
If $f(0)=1$, then $Tf(0)=1$, and the following holds if the integrals exist \begin{equation} \label{eqn:lemma2-1}
\sqrt{p(Z^2)}\int_{\R^2} D\Phi\ e^{-(Z-\Phi)^2}\sqrt{p(\Phi^2)} f(\Phi^2)
= Tf(Z^2). 
\end{equation}

\item[(ii)]
If $f>0$ pointwise, then $Qf>0$ pointwise, and the following holds if the integrals exist \begin{equation}\label{eqn:lemma2-2}
\sqrt{p(Z^2)}\int_{\R^2} D\Phi\ \bar\phi e^{-(Z-\Phi)^2}\sqrt{p(\Phi^2)} f(\Phi^2) 
= \bar \zeta Qf(Z^2). 
\end{equation}
\end{enumerate}
\end{lemma}

\begin{proof}
(i) By definition of $T$ in~\eqref{def:T} and by the Taylor expansion of the kernel~\eqref{eqn:2.6}, 
\begin{align}
Tf(0) = \sqrt{p(0)} \cdot 1 + \int_0^\infty f(s) \cdot 0\ ds = 1. 
\end{align}
By \cite[Proposition 2.5]{BS2020}, \begin{align} \label{eqn:A18}
\int_{\R^2} D\Phi\ e^{-(Z-\Phi)^2}\sqrt{p(\Phi^2)} f(\Phi^2) = e^{-V(Z^2)}, \end{align}
where 
\begin{align}\label{eqn:A19}
e^{-V(t)} &= e^{-t} \( \sqrt{p(0)} f(0) + v(t) \),  \\
v(t) &= \int_0^\infty \sqrt{p(s)} f(s)  e^{-s} I_1(2\sqrt{st})\sqrt{\tfrac{t}{s}}\ ds.
\end{align}
We multiply equation~\eqref{eqn:A19} by $\sqrt{p(t)}$ and use $\sqrt{ p(0)}=f(0)=1$, then \begin{align}
\sqrt{p(t)}e^{-V(t)} 
&= \sqrt{p(t)} e^{-t}+ 
\int_0^\infty f(s) \sqrt{p(t)}\sqrt{p(s)}  e^{-t} e^{-s} I_1(2\sqrt{st})\sqrt{\tfrac{t}{s}}\ ds
 \\
&= Tf(t). \nonumber
\end{align}
Substituting $Z^2$ into $t$ gives the desired (\ref{eqn:lemma2-1}). 

\smallskip \noindent
(ii) Recall $Qf(t) = \int_0^\infty f(s) k_0(t,s) ds$ and $k_0(t,s) > 0$ for all $t$. Since $f>0$ pointwise, we get $Qf(t)>0$ pointwise too.  
By \cite[Lemma 2.6]{BS2020}, \begin{align} \label{eqn:A.7}
\int_{\R^2} D\Phi\ \bar\phi e^{-(Z-\Phi)^2}\sqrt{p(\Phi^2)} f(\Phi^2) 
= \bar \zeta (1- V'(Z^2)) e^{-V(Z^2)},
\end{align} where $V$ is the same as in (\ref{eqn:A19}). 
Since $V(t) = t- \log(f(0) + v(t))$, differentiating gives $1-V'(t) =  \frac{v'(t)}{f(0)+v(t)}$. Hence, using equation~\eqref{eqn:A19} and $\sqrt{p(0)}=1$, 
\begin{align} \label{eqn:A.8}
(1-V'(t))e^{-V(t)} = e^{-t} v'(t)
&= e^{-t} \int_0^\infty \sqrt{p(s)} f(s)  e^{-s}  \frac{\del}{\del t} \( I_1(2\sqrt{st})\sqrt{\tfrac{t}{s}} \) ds  \\
&= e^{-t}   \int_0^\infty \sqrt{p(s)} f(s) e^{-s} I_0(2\sqrt{st})\ ds,  \nonumber
\end{align} 
where the last equality is by Taylor series of the modified Bessel functions \begin{align}
I_0(2\sqrt{st}) = \sum_{m=0}^\infty \frac{s^m t^m}{m! m!}, \qquad
I_1(2\sqrt{st})\sqrt{\tfrac t s} = \sum_{m=0}^\infty \frac{s^m t^{m+1}}{m!( m+1)!}. 
\end{align}
Multiplying equation \eqref{eqn:A.8} by $\sqrt{p(t)}$, we get \begin{align}
\sqrt{p(t)} (1-V'(t))e^{-V(t)} 
&= \int_0^\infty f(s) \sqrt{p(t)}\sqrt{p(s)} e^{-t} e^{-s} I_0(2\sqrt{st})\ ds
= Qf(t). \end{align}
Substituting $Z^2$ into $t$, then plugging into $\sqrt{p(Z^2)} \cdot \eqref{eqn:A.7}$ gives the desired (\ref{eqn:lemma2-2}). 
\end{proof}

\begin{proof}[Proof of Proposition~\ref{prop:Gij^N}]
Let $-N \le i \le j \le N$. 
We first prove the statement for $\nu\in \R$. 
Let $\Del_N$ denote the generator of the random walk on $[-N, N]$, then
\begin{align}
(\Phi, -\Del_N\Phi) &= \sum_{x=-N}^N \( \phi_x(-\Del_N\bar\phi)_x + \psi_x(-\Del_N\bar\psi)_x \) 
= \sum_{x=-N}^{N-1} (\Phi_{x+1} - \Phi_x)^2.
\end{align}
Using assumption~\eqref{A3} and Theorem~\ref{theorem:isomorphism}, 
$G_{ij}^N$ (defined in \eqref{def:Gij^N}) can be expressed as 
\begin{align} \label{eqn:A27}
G_{ij}^N(g, \nu) 
&= \int_0^\infty E_i^N \( \prod_{x=-N}^{N} p(L_{T,x}) \1_{X(T)=j} \) dT  \\
&= \int_{\R^{2(2N+1)}} \bar\phi_i \phi_j e^{-(\Phi, -\Del_N\Phi)} \prod_{x=-N}^{N} p(\Phi_x^2)  \nonumber\\
&=  \int_{\R^{2(2N+1)}} \bar\phi_i \phi_j 
\prod_{x=-N}^{N-1} e^{-(\Phi_{x+1} - \Phi_x)^2}  \prod_{x=-N}^N p(\Phi_x^2) .\nonumber 
\end{align}
We use Lemma~\ref{lemma:outside_integral} to calculate this integral iteratively. First, we decompose $p(\Phi_x^2) = \sqrt{ p(\Phi_x^2)} \sqrt{ p(\Phi_x^2)}$ for all $x$. 
Starting from $-N$, if $i > -N$, we take one of the $\sqrt{ p(\Phi_{-N+1}^2)}$ terms, all $\Phi_{-N}^2$ terms, and calculate the $\Phi_{-N}$ integral. This matches Lemma~\ref{lemma:outside_integral}(i) with $Z = \Phi_{-N+1}$, $\Phi = \Phi_{-N}$, and $f = \sqrt{p}$, giving
\begin{align}
\sqrt{p(\Phi_{-N+1}^2)}\int_{\R^2} D\Phi_{-N}\ e^{ - (\Phi_{-N+1}-\Phi_{-N})^2 } &\sqrt{p(\Phi_{-N}^2)}  \sqrt{p(\Phi_{-N}^2)} 
= T[\sqrt p] (\Phi_{-N+1}^2). 
\end{align}
The process continues until we reach $i$. 
If $j<N$, we also start from $N$ and integrate out $\Phi_N, \dots, \Phi_{j+1}$. 
This gives
\begin{multline}
G_{ij}^N(g, \nu)  =
\int_{\R^{2(j-i+1)}} \bar\phi_i \phi_j 
\cdot T^{N+i}[\sqrt p] (\Phi_{i}^2) \cdot T^{N-j}[\sqrt p] (\Phi_{j}^2) \\
\cdot \Bigg[ \prod_{x=i}^{j-1} e^{-(\Phi_{x+1} - \Phi_x)^2}  \Bigg]
\sqrt{p(\Phi_i^2)} 
\Bigg[ \prod_{x=i+1}^{j-1} p(\Phi_x^2)  \Bigg] \sqrt{p(\Phi_j^2)} .
\end{multline}
We then integrate from $i$ up to $j-1$ using Lemma~\ref{lemma:outside_integral}(ii). The only difference is that there is an extra boson $\bar \phi_i$ that gets carried along. We get 
\begin{align}
G_{ij}^N(g, \nu) 
&=  \int_{\R^2}  \bar \phi_j \phi_j\cdot Q^{j-i}T^{N+i}[\sqrt p] (\Phi_{j}^2)\cdot  T^{N-j}[\sqrt p] (\Phi_{j}^2). 
\end{align}
By Theorem~\ref{theorem:supersymmetry} and the exponential decay of $\sqrt{p}$, \begin{align} \label{eqn:25}
G_{ij}^N(g, \nu) 
&=  \int_{\R^2}  \bar \psi_j \psi_j \cdot Q^{j-i}T^{N+i}[\sqrt p] (\Phi_{j}^2)\cdot  T^{N-j}[\sqrt p] (\Phi_{j}^2)  \\
&=  \int_{\R^2} Q^{j-i}T^{N+i}[\sqrt p] (u^2 + v^2)\cdot  T^{N-j}[\sqrt p] (u^2+v^2)
\frac 1 \pi du dv\nonumber  \\
&= \int_0^\infty Q^{j-i}T^{N+i}[\sqrt p] (t)\cdot  T^{N-j}[\sqrt p] (t) dt \nonumber \\
&= \left \langle Q^{j-i}T^{N+i}[\sqrt p] , \overline{ \Big( T^{N-j}[\sqrt p] \Big)}\right \rangle,\nonumber
\end{align}
as desired. 

For complex $\nu$, observe that both sides of (\ref{eqn:25}) are defined and holomorphic in $\nu\in\C$. We get the result by the uniqueness of analytic continuation. 
\end{proof}

\section{Proof of Lemma~\ref{lemma:Hn}} \label{appendix:Hn}
The proof of Lemma~\ref{lemma:Hn} is via a direct computation. 
\begin{proof}
We first calculate derivatives of $\lambda(g, \nu) = \norm{Q(g, \nu)}$, by viewing $Q(g, \nu)$ as a perturbation of $Q(g_0, \nu_0)$. This calculation is similar to that of the Rayleigh--Schr\"odinger series. 
Then we calculate the derivatives of the function $H_n$ defined in \eqref{def:Hn}, just using differentiation rules. 
When the derivatives are evaluated at $(g_0, \nu_0)$, the formulas will simplify because  $Q(g_0, \nu_0)h_0 = \lambda(g_0, \nu_0) h_0$, and the claimed results will become apparent.

Since $\lambda(g, \nu) = \norm{Q(g, \nu)}$ is an isolated simple eigenvalue (Lemma~\ref{lemma:Q}), there exists $\delta>0$ such that $\lambda_0 = \lambda(g_0, \nu_0)$ is distance $2\delta$ away from the rest of the spectrum of $Q(g_0, \nu_0)$. 
For $(g, \nu)$ near $(g_0, \nu_0)$, 
the projection operator to the eigenspace $E_{\lambda(g, \nu)}$ of $Q(g, \nu)$ is given by
\begin{align}
P(g, \nu) = -\frac 1 {2\pi i} \oint_{\abs{\lambda_0 - \zeta}=\delta} (Q(g, \nu) - \zeta)\inv d\zeta .
\end{align}
Thus, we have $QPh_0 = \lambda Ph_0$, and 
\begin{align}\label{eqn:B2}
\lambda=  \frac{ \langle QPh_0, h_0 \rangle} {\langle Ph_0, h_0 \rangle}. 
\end{align}
We differentiate this equation. Note the dependence on $g$ and $\nu$ only come from $Q$ and $P$. 

We use subscripts to denote partial derivatives. The $\nu$-derivative of \eqref{eqn:B2} is \begin{align} \label{eqn:B3}
\lambda_\nu =
\frac{ \langle Q_\nu Ph_0, h_0 \rangle + \langle QP_\nu h_0, h_0 \rangle}{\langle Ph_0, h_0 \rangle}
-\frac{ \langle QPh_0, h_0 \rangle \langle P_\nu h_0, h_0 \rangle}{\langle Ph_0, h_0 \rangle^2}. 
\end{align}
When evaluated at $(g_0, \nu_0)$, we know $Ph_0=h_0$, $Qh_0 =\lambda_0 h_0$, $\norm{h_0}_2 = 1$, and $Q$ is self-adjoint. Hence, 
\begin{align}
\lambda_\nu \rvert_{g_0, \nu_0}
= \langle Q_\nu h_0, h_0 \rangle + \langle P_\nu h_0, Q h_0 \rangle
- \lambda_0 \langle P_\nu h_0, h_0 \rangle
= \langle Q_\nu h_0, h_0 \rangle. 
\end{align}
Similarly, $\lambda_g  \rvert_{g_0, \nu_0} = \langle Q_g h_0, h_0 \rangle$. 

For second derivatives, we let $*=g, \nu$, differentiate \eqref{eqn:B3}, and then evaluate at $(g_0, \nu_0)$. This gives \begin{align} \label{eqn:B5}
\lambda_{\nu  *} \rvert_{g_0, \nu_0} = 
&\ \langle Q_{\nu *} h_0, h_0 \rangle 
+ \langle Q_\nu P_* h_0, h_0 \rangle + \langle Q_* P_\nu h_0, h_0 \rangle  \\
&- \langle Q_\nu h_0, h_0 \rangle \langle P_* h_0, h_0 \rangle - \langle Q_* h_0, h_0 \rangle \langle P_\nu h_0, h_0 \rangle. \nonumber
\end{align}
We claim $\langle P_* h_0, h_0 \rangle = 0$. This is because \begin{align} \label{eqn:B6}
\langle P_* h_0, h_0 \rangle
&= \frac{1}{2\pi i}  \oint_{\abs{\lambda_0 - \zeta}=\delta} \langle (Q-\zeta)\inv Q_* (Q-\zeta)\inv h_0, h_0  \rangle d\zeta \\
&= \frac{1}{2\pi i}  \oint_{\abs{\lambda_0 - \zeta}=\delta} 
\frac 1 {\lambda_0 - \zeta }
\langle  Q_*  h_0, (Q-\bar \zeta)\inv h_0  \rangle d\zeta \nonumber\\
&= \frac{1}{2\pi i} \langle  Q_*  h_0, h_0  \rangle \oint_{\abs{\lambda_0 - \zeta}=\delta} 
\frac 1 { ( \lambda_0 - \zeta )^2 }
 d\zeta  = 0.\nonumber
\end{align}
Next, we calculate $\langle Q_\nu P_* h_0, h_0 \rangle$ in equation~\eqref{eqn:B5}. Since $Q$ is self-adjoint, by the spectral theorem, there exists a real orthonormal eigenbasis $\{(\mu_j, \psi_j)\}_j$ of $Q(g_0, \nu_0)$. Using these, we decompose \begin{align}
Q_*h_0 = \sum_j \langle Q_* h_0, \psi_j \rangle \psi_j,
\end{align}
so \begin{align} \label{eqn:B8}
\langle Q_\nu P_* h_0, h_0 \rangle
&= \frac{1}{2\pi i}  \oint_{\abs{\lambda_0 - \zeta}=\delta} \langle Q_\nu (Q-\zeta)\inv Q_* (Q-\zeta)\inv h_0, h_0  \rangle d\zeta \\
&= \frac{1}{2\pi i}  \oint_{\abs{\lambda_0 - \zeta}=\delta} \frac{1}{\lambda_0 - \zeta}
\langle  (Q-\zeta)\inv Q_*  h_0, Q_\nu h_0  \rangle d\zeta \nonumber\\
&= \frac{1}{2\pi i}  \oint_{\abs{\lambda_0 - \zeta}=\delta} \frac{1}{\lambda_0 - \zeta} 
\sum_j \frac{ \langle Q_* h_0, \psi_j \rangle \langle Q_\nu h_0, \psi_j \rangle }{ \mu_j - \zeta} d\zeta \nonumber\\
&= - \sum_{\mu_j \ne \lambda_0} \frac{ \langle Q_* h_0, \psi_j \rangle \langle Q_\nu h_0, \psi_j \rangle }{ \mu_j - \lambda_0}.  \nonumber
\end{align}
In the last equality, the $\mu_j = \lambda_0$ term vanishes for the same reason as in~\eqref{eqn:B6}. To write this more compactly, we define $P^\perp = I - P$ where $I$ is the identity operator, then equation~\eqref{eqn:B8} can be written as \begin{align}
\langle Q_\nu P_* h_0, h_0 \rangle
= \langle (\lambda_0 - Q)\inv P^\perp Q_* h_0, P^\perp Q_\nu h_0 \rangle. 
\end{align}
We also have $
\langle Q_* P_\nu h_0, h_0 \rangle$ equal to the same expression, by the symmetry between $*$ and $\nu$ in equation~\eqref{eqn:B8}. 
Putting together, equation~\eqref{eqn:B5} simplifies to \begin{align}
\lambda_{\nu  *} \rvert_{g_0, \nu_0} = 
&\ \langle Q_{\nu *} h_0, h_0 \rangle 
+2 \langle (\lambda_0 - Q)\inv P^\perp Q_* h_0, P^\perp Q_\nu h_0 \rangle. \label{eqn:B10}
\end{align}

We next turn to the derivatives of $H_n(g, \nu) = \langle Q^n(g, \nu)h_0, h_0 \rangle^{1/n}$. A direct computation gives 
\begin{align}
\del_\nu H_{n} 
&=\frac{1}{n} \langle Q^n h_0, h_0 \rangle^{1/n-1}   
 \langle (Q^n)_\nu h_0, h_0 \rangle, \\
\del_g\del_\nu H_{n} 
&= \frac{1}{n} \(\frac{1}{n} - 1\) \langle Q^n h_0, h_0 \rangle^{1/n-2}  
\langle (Q^n)_g h_0, h_0 \rangle 
\langle (Q^n)_\nu h_0, h_0 \rangle \\
&\qquad+ \frac{1}{n} \langle Q^n h_0, h_0 \rangle^{1/n-1}   
\langle (Q^n)_{\nu g} h_0, h_0 \rangle. \nonumber
\end{align} 
When evaluated at $(g_0, \nu_0)$, we have $Q(g_0, \nu_0)h_0 = \lambda_0h_0$ with $\lambda_0 =1$. This gives 
\begin{align}
\del_\nu H_{n} (g_0, \nu_0)
&= \frac{1}{n}  \langle (Q^n)_\nu h_0, h_0 \rangle 
= \langle Q_\nu h_0, h_0 \rangle = \lambda_\nu, \\
\del_g\del_\nu H_{n} (g_0, \nu_0) 
&= \frac{1}{n} \(\frac{1}{n} - 1\) \langle (Q^n)_g h_0, h_0 \rangle 
\langle (Q^n)_\nu h_0, h_0 \rangle + \frac{1}{n} 
 \langle (Q^n)_{\nu g} h_0, h_0 \rangle   \\
&= (1-n) \langle Q_g h_0, h_0 \rangle 
\langle Q_\nu h_0, h_0 \rangle 
+  \langle Q_{\nu g} h_0, h_0 \rangle \nonumber\\
&\qquad +\frac{2}{n} \sum_{0 \le i< j \le n} \langle Q^{i-1}Q_\nu Q^{j-i-1} Q_g Q^{n-j} h_0, h_0 \rangle \nonumber \\
&=(1-n) \langle Q_g h_0, h_0 \rangle 
\langle Q_\nu h_0, h_0 \rangle 
+  \langle Q_{\nu g} h_0, h_0 \rangle \nonumber\\
&\qquad+\frac{2}{n} \sum_{0 \le i< j \le n} \langle Q^{j-i-1} Q_g h_0,Q_\nu h_0 \rangle. \nonumber
\end{align} 
For the last sum, we decompose $Q_g h_0 = PQ_g h_0 + P^\perp Q_g h_0$ and similarly for $Q_\nu h_0$. The parts that are in the eigenspace $E_1$ sum to cancel with $(1-n) \langle Q_g h_0, h_0 \rangle \langle Q_\nu h_0, h_0 \rangle $ exactly, leaving \begin{align}
\del_g\del_\nu H_{n} (g_0, \nu_0) =
\langle Q_{\nu g} h_0, h_0 \rangle 
+\frac{2}{n} \sum_{0 \le i< j \le n} \langle Q^{j-i-1} P^\perp Q_g h_0, P^\perp Q_\nu h_0 \rangle. 
\end{align}
Summing diagonally, as $n\to\infty$ we get 
\begin{align}
\del_g\del_\nu H_{n} (g_0, \nu_0) &\to
\langle Q_{\nu g} h_0, h_0 \rangle + 2 \langle (1 - Q)\inv P^\perp Q_g h_0, P^\perp Q_\nu h_0 \rangle
= \lambda_{\nu g}\rvert_{g_0, \nu_0} ,
\end{align}
by the first computation \eqref{eqn:B10}. 
The calculation for $\lambda_{\nu\nu}\rvert_{g_0, \nu_0} $ is analogous. 
\end{proof}

\end{appendices}

\bibliographystyle{plain}

\end{document}